\long\def\history#1{}
\history{
Original submitted, Dec. 96.
Two changes suggested by I. Hambleton, Jan. 97. 
Suggestions by Teichner Fall 97.
Referee's report March 1998
}
\def\START{}
\magnification1200

\gdef\namedis#1{}

\def\splitthm#1\ #2@{\gdef\ThmKind{#1}\gdef\ThmNum{#2}}
\def\twolabel#1#2{\edef\xx{#1@}\expandafter\splitthm\xx\relax%
\let\ThmKindA=\ThmKind\let\ThmNumA=\ThmNum%
\edef\xx{#2@}\expandafter\splitthm\xx\relax%
\ThmKind s \ThmNumA\ and \ThmNum}

%
\def\topmatter{\bgroup\bgroup\baselineskip=25pt}
\def\endtopmatter{\null\egroup\egroup}
\def\title#1{\centerline{\bf#1}}
\def\authors#1{\centerline{\bf#1}}
\def\dedication#1{\centerline{\sl#1}}
\def\thanks#1{\footnote{}{#1}}
\newbox\endbox
\global\let\hhh=\relax
\global\let\HHH=\hfill
\newdimen\dom
\setbox0=\hbox{Department of Mathematics}
\global\dom=\wd0
\def\address#1{{\parindent=0pt\parskip=0pt\baselineskip=8pt%
\def\\{\par}\global\setbox\endbox=\hbox to\hsize{%
\unhbox\endbox\hhh\hbox to \dom{\vtop{\hsize=2.5in#1}\hss}\HHH}%
\global\let\hhh=\hfill\global\let\HHH=\relax}}
\let\lrtend=\end
\def\end{\vskip 10pt\noindent \box\endbox\lrtend}
\def\evenheader{}
\def\oddheader{}
\font\SmallCaps=cmr9
\font\VerySmallCaps=cmr7
\def\upper#1#2@{\SmallCaps #1\VerySmallCaps\uppercase{#2}}
\headline{\ifnum\pageno=1\relax\hfill\else%
{\hfill\ifodd\pageno\oddheader\else\evenheader\fi\hfill}\fi}

\newwrite\labellist
\newif\ifwritelabels\writelabelsfalse
\def\makelabellist{\immediate\openout\labellist=\jobname.labelsA
\writelabelstrue}

\font\Bbb=msbm10
\def\bbb#1{\hbox{\Bbb #1}}

\newdimen\bWid\bWid=.9in
\newdimen\HSIZE
\newdimen\dEP
\newdimen\HAL

\def\spin{$Spin$}
\def\semidirect{\hbox{%
$\vrule width .3pt height 4.5pt depth 0pt\hskip-1.75pt\times$}}

\def\WHATSIT#1\ #2@{\gdef\first{#1}\gdef\second{#2}}
\def\whatsit#1{\setbox0=\hbox{\expandafter\WHATSIT#1@}}
\def\twotheorems#1#2{%
\whatsit#1%
\let\firstA=\first\let\secondA=\second%
\whatsit#2%
\ifx\firstA\first\relax\first s\ \secondA\ and \second\else%
#1\ and #2\fi
}

\def\hmtytwo{\hmty_{0}}
\def\hmtybig{\hmty_{\%}}

\def\normact{\bullet}
\def\autact{\bullet}
\def\KSact{\bullet}
\def\hmtyone{\hmty_{1}}
\def\and{\ \&\ }
\def\gtop{G/TOP}
\def\gdiff{G/O}
\def\gcat{G/\CAT}
\def\topdiff{TOP/O}
\def\Ncat{{\cal L}^{\CAT}}
\def\Ntop{{\cal L}^{\TOP}}

\def\wtNcat#1{\hbox{$^{#1}$}\hskip-1pt\widetilde\Ncat}

\def\Scat{{\cal S}^{CAT}}
\def\Stop{{\cal S}^{{\TOP}}}
\def\Sdiff{{\cal S}^{{\DIFF}}}
\def\barScat{\bar\Scat}
\def\barStop{\bar\Stop}
\def\barSdiff{\bar\Sdiff}
\def\wtScat#1{\hbox{$^{#1}$}\hskip-1pt\widetilde\Scat}

\def\wtSdiff#1{\hbox{$^{#1}$}\hskip-1pt\widetilde\Sdiff}
\def\WSEcat{W\hskip-2ptSE^{\CAT}}
\def\WSEtop{W\hskip-2ptSE^{\TOP}}
\def\WSEdiff{W\hskip-2ptSE^{\DIFF}}

\def\remark{\noindent{\bf Remark:}\ }
\def\remarks{\noindent{\bf Remarks:}\ }

\def\disjointunion{\perp\hskip -2pt\perp}
\def\Z{{\bf Z}}
\def\Q{{\bf Q}}
\def\R{{\bf R}}
\def\cy#1{\Z/{#1}\Z}

\def\rel{{\rm rel}}
\def\CAT{CAT}
\def\TOP{TOP}
\def\DIFF{DIFF}
\def\PL{PL}
\def\comp{\circ}

\def\N{N}

\def\image{{\rm Im}}
\def\hmty{H\hskip-1pt E^{^+}}
\def\CP{C\hskip -1pt P}
\def\chern{\widehat{\CP}^2}

\def\RA#1{\setbox0=\hbox{$\scriptstyle11#1$}%
\displaystyle\mathop{\ \hbox to \wd0 {\rightarrowfill}}^{#1}\ %
}
\def\rzb#1{\hbox to 0pt{$\scriptstyle#1$\hss}}
\def\lzb#1{\hbox to 0pt{\hss$\scriptstyle#1$}}
\def\LWG#1#2{{L}^{^{#1}}_{#2}}
\def\overo#1{\mathchoice%
{\displaystyle{\mathop{#1}^{\circ}_{\vphantom{\circ}}}}
{\displaystyle{\mathop{\hbox{$#1$}}^{\circ}_{\vphantom{\circ}}}}
{\displaystyle{\mathop{\scriptstyle#1}^{\scriptscriptstyle\circ}_{\vphantom{\scriptscriptstyle\circ}}}}
{}%
}
\def\br#1#2{[\,{#1},{#2}\,]}
\def\lct{{\bf L}\langle1\rangle}
\def\zpi#1#2{\ifnum#1=1 \relax\Z{#2}\else\Z[{#2}]\fi}

\def\Ndef#1#2{%
\gdef\next{#2}%
\ifx\next\empty\xdef#1{\number\displaynum\namedis#1}%
{{\ifwritelabels\immediate\write\labellist%
{\def\string #1{\number\displaynum}}\fi}}%
\else%
\xdef#1{#2\ \number\displaynum\namedis#1}%
{{\ifwritelabels\immediate\write\labellist%
{\def\string #1{#2\ \number\displaynum}}\fi}}%
\fi%
\global\advance\displaynum by1}

\newcount\displaynum\global\displaynum=1
\def\Ddef#1#2{\xdef#1{#2\namedis#1}}

\newcount\sectionnumber\global\sectionnumber=0

\def\remembersection#1{%
{\ifwritelabels\immediate\write\labellist%
{\def\string #1{\number\sectionnumber}}\fi}
\gdef#1{\number\sectionnumber}%
}
\def\Remembersection#1{%
\global\advance\sectionnumber by -1 %
{\ifwritelabels\immediate\write\labellist%
{\def\string #1{\number\sectionnumber. \xx}}\fi}
\global\advance\sectionnumber by 1 %
}
\newcount\truesection
\def\newsec#1{%
\medskip\penalty-100\noindent{\bf\S\number\sectionnumber. #1}%
\xdef\xx{#1}%
\global\truesection=\sectionnumber
\medskip\global\advance\sectionnumber by 1}
\catcode`?=11
\def\stringit#1#2@{\def\xx{#2}%
\edef\rr{\expandafter\noexpand\csname ?\xx\endcsname}}
\def\cite#1{[#1]\expandafter\stringit\string#1@\removeequivalent#1\rr}
\def\ecite#1#2{[#1, #2\hskip 2pt]%
\expandafter\stringit\string#1@\removeequivalent#1\rr}
\let\?=\relax
\def\biblist{}
\toksdef\t?a=0 \toksdef\t?b=2 %
\def\append#1{\t?a={\\#1}\t?b=\expandafter{\biblist}%
\edef\biblist{\the\t?b\the\t?a}}
\def\?biblist{}
\def\appendA#1@{%
\t?a=\expandafter{\csname ??#1\endcsname}\t?b=\expandafter{\?biblist}%
\ifx\?biblist\empty\edef\?biblist{\?\the\t?a}\else
\edef\?biblist{\the\t?b\?\the\t?a}\fi}%
\def\removeequivalent#1#2{\let\given=#1\def\givenb{#2}%
\ifx\biblist\empty\else\edef\biblist{\expandafter\plugh\biblist\plugh}\fi}
\def\plugh\\#1#2{\ifx#1\given%
\noexpand\\\givenb%
\else\noexpand\\\noexpand#1\fi%
\ifx#2\plugh\hgulp\fi\plugh#2}
\def\hgulp\fi\plugh\plugh{\fi}
\def\showbiblist{\ifx\biblist\empty\else\let\\=\sx?x\relax%
\biblist\let\?=\f?orm\?biblist\fi}

\newdimen\fixsp
\def\f?orm#1{\edef\xx{#1?}\expandafter\f??orm\xx}
\def\f??orm[#1]#2?{\clear?em\form[#1]#2!}
\def\stringitA#1#2#3@{\def\xx{#2}
\ifx\xx\qm\relax\appendA#3@\fi}%
\def\sx?x#1{\expandafter\stringitA\string#1@}

\newcount\bcount\bcount=1

\def\stringitB#1#2@#3@{\def\xx{#2}%
\expandafter\xdef\csname ??\xx\endcsname{#3}}

\fixsp=.4em
\def\same#1{%
\setbox0=\hbox{\spaceskip=\fixsp#1}\authdim=\wd0 %
\global\authtoks{{\advance\authdim by -2pt%
\vrule width\authdim height 0pt depth .1pt\hskip2pt}}
}%
\def\sameand#1#2{%
\setbox0=\hbox{\spaceskip=\fixsp#1}\authdim=\wd0 %
\global\authtoks{{\advance\authdim by -4pt%
\vrule width\authdim height 0pt depth .1pt\hskip4pt#2}}
}%
\def\Bdef#1#2#3{\edef#1{\number\bcount\namedis{#1}}%
\append#1%
\expandafter\stringitB\string #1@[#1]{#3}@%
\global\advance\bcount by1}

\def\clear?em{%
\global\authtoks{}%
\global\pubtoks{}%
\global\articletoks{}%
\global\booktoks{}%
\global\inbooktoks{}%
\global\journaltoks{}%
\global\voltoks{}%
\global\editortoks{}%
\global\pagestoks{}%
\global\yeartoks{}%
\global\extratoks{}%
\global\preprinttoks{}%
}

\catcode`?=12
\def\qm{?}
\catcode`?=11
\newdimen\authdim
\newdimen\boff\newdimen\bsize
\boff=20pt\bsize=\hsize\advance\bsize by-\boff
\newbox\authbox
\def\form[#1]#2!{\par\noindent\hbox to\hsize{%
\vtop{\hsize=\boff\parindent=0pt\hbox to\boff{\hfil[#1]\hfil}\vfil}%
\vtop{\parindent=0pt\hsize=\bsize%
#2\displayref%
\vfil}\hfil}\par\vfil%
}%
\catcode`?=12
\newskip\commaskip\commaskip=.33em plus 12pt minus 2pt
\newcount\reftype
\def\displayref{{{\spaceskip=\fixsp\the\authtoks}},\hskip\commaskip
\edef\xx{\the\pubtoks}\ifx\xx\empty\else%
\global\expandafter\pubtoks{\xx, }\fi%
\edef\xy{\the\editortoks}\ifx\xy\empty\else%
\global\expandafter\editortoks{\xy, }\fi%
\edef\xz{\the\extratoks}\ifx\xz\empty\else%
\global\expandafter\extratoks{\xz, }\fi%
\ifcase\reftype%
\the\articletoks. {\sl\the\journaltoks}\ {\bf\the\voltoks} (\the\yeartoks), \the\pagestoks
\or{\sl\the\booktoks}, \the\extratoks\the\pubtoks \the\yeartoks%
\or\the\articletoks. In \the\editortoks {\sl \the\inbooktoks}, \the\extratoks
\edef\xt{\the\pagestoks}\ifx\xt\empty\else pp. \the\pagestoks, \fi
\the\pubtoks\the\yeartoks
\or\the\articletoks, \the\preprinttoks
\else\fi.\vskip4pt}

\newtoks\authtoks
\newtoks\pubtoks
\newtoks\articletoks
\newtoks\booktoks
\newtoks\inbooktoks
\newtoks\journaltoks
\newtoks\voltoks
\newtoks\editortoks
\newtoks\pagestoks
\newtoks\yeartoks
\newtoks\extratoks
\newtoks\preprinttoks

\def\article#1{\global\articletoks{#1}}
\def\pub#1{\global\pubtoks{#1}}
\def\editor#1{\global\editortoks{#1}}
\def\book#1{\global\booktoks{#1}\global\reftype=1}
\def\preprint#1{\global\preprinttoks{#1}\global\reftype=3}
\def\inbook#1{\global\inbooktoks{#1}\global\reftype=2}
\def\journal#1{\global\journaltoks{#1}\global\reftype=0}
\def\vol#1{\global\voltoks{#1}}
\def\pages#1-#2.{\global\pagestoks{#1{--}#2}}
\def\author#1{\global\authtoks{#1}}
\def\yr#1{\global\yeartoks{#1}}
\def\extra#1{\global\extratoks{#1}}


\Bdef\bAK{}{\author{S.\ Akbulut}
\yr{1984}
\article{A fake 4{--}manifold}
\inbook{Four{--}Manifold Theory (Durham, NH, 1982)}
\editor{C.\ McA.\ Gordon and R.\ C.\ Kirby (Eds.)}
\extra{Volume 35 of {\sl Contemp.\ Math.}}
\pub{Amer.\ Math.\ Soc.}
\pages 75-142.
}
\Bdef\bACork{}{\author{S.\ Akbulut}
\same{S.\ Akbulut}
\article{A fake compact contractible 4{--}manifold}
\journal{J.\ Differential Geom.}
\vol{33}
\yr{1991}
\pages 335-356.
}
\Bdef\bAM{AM}{ \author{S.\ Akbulut and J.\ D.\ McCarthy}
\sameand{S.\ Akbulut}{and J.\ D.\ McCarthy}
\yr{1990}
\book{Casson's Invariant for Oriented Homology 3{--}Spheres}
\extra{Volume 36 of {\sl Math.\ Notes.}}
\pub{Princeton Univ.\ Press}}
\Bdef\bAsk{AM}{ \author{N.\ Askitas}
\article{Embeddings of $2$-spheres in $4$-manifolds}
\journal{Manuscripta Math.}
\yr{1996}
\vol{89}
\pages 35-47.}
\Bdef\bAskb{AM}{ \author{N.\ Askitas}
\article{Erratum: Embeddings of $2$-spheres in $4$-manifolds}
\journal{Manuscripta Math.}
\yr{1996}
\vol{90}
\pages 137-138.}
\Bdef\bBV{BV}{\author{J.\ M.\ Boardman and R.\ M.\ Vogt}
\yr{1973}
\book{Homotopy Invariant Algebraic Structures on Topological Spaces}
\extra{Volume 347 of {\sl Lecture Notes in Math.}}
\pub{Springer{--}Verlag}}
\Bdef\bBT{BT}{\author{W.\ Browder}
\yr{1961}
\article{Torsion in $H${--}spaces}
\journal{Ann.\ of Math.}
\vol{74}
\pages 24-51.}
\Bdef\bBB{B}{\author{W.\ Browder}
\same{W.\ Browder}
\yr{1972}
\book{Surgery on Simply Connected Manifolds}
\extra{Ergib.\ der Math.\ 65}
\pub{Springer{--}Verlag}}
\Bdef\bGBr{CSs}{\author{G.\ Brumfiel}
\yr{1971}
\article{Homotopy equivalences of almost smooth manifolds}
\journal{Comment.\ Math.\ Helv.\ }
\vol{46}
\pages 381-407.}%
\Bdef\bCSs{CSs}{\author{S.\ E.\ Cappell and J.\ L.\ Shaneson}
\yr{1971}
\article{On four dimensional surgery and applications}
\journal{Comment.\ Math.\ Helv.\ }
\vol{46}
\pages 500-528.}

\Bdef\bCSz{CSz}{\author{S.\ E.\ Cappell and J.\ L.\ Shaneson}
\same{S.\ E.\ Cappell and J.\ L.\ Shaneson}
\yr{1974}
\article{The codimension two placement problem and homology equivalent
manifolds}
\journal{Ann.\ of Math.}
\vol{99}
\pages 277-348.}
\Bdef\bCSfh{CSh}{\author{S.\ E.\ Cappell and J.\ L.\ Shaneson}
\same{S.\ E.\ Cappell and J.\ L.\ Shaneson}
\yr{1976}
\article{Some new four{--}manifolds}
\journal{Ann.\ of Math.}
\vol{104}
\pages 61-72.}%
\Bdef\bCSh{CSh}{\author{S.\ E.\ Cappell and J.\ L.\ Shaneson}
\same{S.\ E.\ Cappell and J.\ L.\ Shaneson}
\yr{1985}
\article{On $4${--}dimensional $s${--}cobordisms}
\journal{J.\ Differential Geom.}
\vol{22}
\pages 97-115. }
\Bdef\bAC{C}{\author{A.\ J.\ Casson }
\yr{1986}
\article{Three lectures on new infinite constructions in 4{--}dimensional
manifolds}
\editor{L.\ Guillou and A.\ Marin (Eds.)}
\inbook{A la Recherche de la Topologie Purdue}
\extra{Volume 62 of {\sl Prog.\ Math.}}
\pub{ Birkh\"auser}
\pages 201-244.}
\Bdef\bCG{C}{\author{A.\ J.\ Casson and C.\ McA.\ Gordon}
\article{On slice knots in dimension three}
\editor{R.\ J.\ Milgram}
\inbook{Algebraic and Geometric Topology (Stanford, 1976)}
\extra{volume 32, part 2 of {\sl Proc.\ Sympos.\ Pure Math.}}
\pub{Amer.\ Math.\ Soc.}
\yr{1978}
\pages 39-53.}
\Bdef\bCGa{C}{\author{A.\ J.\ Casson and C.\ McA.\ Gordon}
\same{A.\ J.\ Casson and C.\ McA.\ Gordon}
\yr{1986}
\article{Cobordism of classical knots}
\editor{L.\ Guillou and A.\ Marin (Eds.)}
\inbook{A la Recherche de la Topologie Purdue}
\extra{Volume 62 of {\sl Prog.\ Math.}}
\pub{ Birkh\"auser}
\pages 181-200.}
\Bdef\bCF{C}{\author{A.\ J.\ Casson and M.\ H.\ Freedman}
\article{Atomic surgery problems}
\inbook{Four-Manifold Theory}
\editor{C.\ McA.\ Gordon and R.\ C.\ Kirby}
\extra{Volume 35 of {\sl Contemp.\ Math.}}
\pub{Amer.\ Math.\ Soc.}
\yr{1984}
\pages 181-199.}
\Bdef\bCf{C}{\author{J.\ Cerf}
\yr{1968}
\book{Sur les diff\'eomorphismes 
de la sph\`ere de dimension trois}
\extra{Volume 53 of {\sl Lecture Notes in Math.}}
\pub{Springer{--}Verlag}}
\Bdef\bCH{CH}{\author{T.\ D.\ Cochran and N.\ Habegger}
\yr{1990}
\article{On the homotopy theory of simply connected four manifolds}
\journal{Topology}
\vol{29}
\pages 419-440.}
\Bdef\bCFHS{}{\author{C.\ L.\ Curtis, M.\ H.\ Freedman, W.-C.\ Hsiang and
R.\ Stong}
\article{A decomposition theorem for $h${--}cobordant 
smooth simply{--}connected compact 4{--}manifolds}
\journal{Invent.\ Math.}
\vol{123}
\yr{1996}
\pages 343-348.
}
\Bdef\bDF{DF}{\author{S.\ De Michelis and M.\ H.\ Freedman}
\yr{1992}
\article{Uncountably many exotic ${\bf R}^4$'s in standard $4${--}space}
\journal{J.\ Differential Geom.}
\pages 219-254.
}
\Bdef\bSDa{D1}{\author{S.\ K.\ Donaldson}
\yr{1983}
\article{An application of gauge theory to 4{--}dimensional topology}
\journal{J.\ Differential Geom.}
\vol{18}
\pages 279-315.}
\Bdef\bSDb{D2}{\author{S.\ K.\ Donaldson}
\same{S.\ K.\ Donaldson}
\yr{1986}
\article{Connections, cohomology and the intersection
 forms of 4{--}manifolds}
\journal{J.\ Differential Geom.}
\vol{24}
\pages 275-341.}
\Bdef\bSDc{D3}{\author{S.\ K.\ Donaldson}
\same{S.\ K.\ Donaldson}
\yr{1987}
\article{Irrationality and the $h${--}cobordism conjecture}
\journal{J.\ Differential Geom.}
\vol{26}
\pages141-168.}
\Bdef\bDK{DK}{\author{S.\ K.\ Donaldson and P.\ B.\ Kronheimer}
\sameand{S.\ K.\ Donaldson}{and P.\ B.\ Kronheimer}
\book{Geometry of Four{--}Manifolds}
\yr{1990}
\pub{Oxford Univ.\ Press}
}
\Bdef\bEM{266}{\author{B.\ Eckmann and P.\ Linnell}
\yr{1983}
\article{Poincar\'e duality groups of dimension 2, II}
\journal{Comment. Math. Helv.}
\vol{58}
\pages 111-114.}
\Bdef\bEL{267}{\author{B.\ Eckmann and H.\ Muller}
\yr{1980}
\article{Poincar\'e duality groups of dimension 2}
\journal{Comment.\ Math. Helv.\ }
\vol{55}
\pages 510-520.}
\Bdef\bFRR{}{\author{S.\ C.\ Ferry, A.\ A.\ Ranicki
and J.\ M.\ Rosenberg, (Editors)}
\book{Novikov Conjectures, Index Theorems and Rigidity,
Volumes 1 and 2, (Oberwolfach 1993)}
\extra{Volume 226 \& 227 of {\sl London Math.\ Soc.\ Lect.\ Note Ser.}}
\pub{Cambridge Univ.\ Press}
\yr{1995}
}
\Bdef\bFS{FS}{\author{R.\ Fintushel and R.\ J.\ Stern}
\article{Knots. links and 4{--}manifolds}
\preprint{Preprint, 1996}
}
\Bdef\bFa{F}{\author{M.\ H.\ Freedman}
\yr{1982}
\article{The topology of four{--}dimensional manifolds}
\journal{J.\ Differential Geom.\ }
\vol{17}
\pages 357-453.}
\Bdef\bFb{F2}{\author{M.\ H.\ Freedman}
\same{M.\ H.\ Freedman}
\yr{1984}
\article{The disk theorem for four dimensional manifolds}
\inbook{Proc.\ International Congress of Mathematicians, Warsaw 1983}
\pages 647-663.
\pub{Polish Scientific Publishers, Warsaw}}
\Bdef\bFQ{FQ}{\author{M.\ H.\ Freedman and F.\ S.\ Quinn}
\sameand{M.\ H.\ Freedman}{and F.\ S.\ Quinn}
\yr{1990}
\book{Topology of 4{--}Manifolds}
\pub{Princeton Univ.\ Press}}
\Bdef\bFTx{FT}{\author{M.\ H.\ Freedman and L.\ R.\ Taylor}
\sameand{M.\ H.\ Freedman}{and L.\ R.\ Taylor}
\article{$\Lambda${--}Splitting 4{--}manifolds}
\journal{Topology}
\vol{16}
\yr{1977}
\pages 181-184.
}
\Bdef\bFT{FT}{\author{M.\ H.\ Freedman and P.\ Teichner}
\sameand{M.\ H.\ Freedman}{and P.\ Teichner}
\yr{1995}
\article{4{--}manifold topology I: Subexponential groups}
\journal{Invent.\ Math.\ }
\vol{122}
\pages 509-529.}
\Bdef\bFM{FM}{\author{R.\ D.\ Friedman and J.\ W.\ Morgan}
\book{Smooth Four{--}Manifolds and Complex Surfaces}
\extra{Volume 27 of {\sl Ergeb.\ Math.\ Grenzgeb.\ (3)}}
\pub{Springer{--}Verlag}
\yr{1994}
}
\Bdef\bFru{F}{\author{M.\ Furuta}
\yr{1995)}
\article{Monopole equation and the 11/8{--}conjecture}
\preprint{Preprint RIMS, Kyoto}}
\Bdef\bGS{GS}{\author{S.\ Gitler and J.\ Stasheff}
\yr{1965}
\article{The first exotic class of $BG$}
\journal{Topology}
\vol{4}
\pages 257-266.}
\Bdef\bBG{G}{\author{R.\ E.\ Gompf}
\yr{1985}
\article{An infinite set of exotic $R^4$'s}
\journal{J.\ Differential Geom.\ }
\vol{21}
\pages 283-300.}
\Bdef\bBGx{G}{\author{R.\ E.\ Gompf}
\same{R.\ E.\ Gompf}
\yr{1991}
\article{Sums of elliptic surfaces}
\journal{J.\ Differential Geom.\ }
\vol{34}
\pages 93-114.}%
\Bdef\bHM{H-M}{\author{I.\ Hambleton and R.\ J.\ Milgram}
\yr{1978}
\article{Poincar\'e transversality for double covers}
\journal{Canad.\ J.\ Math.}
\vol{6}
\pages 1319-1330.}
\Bdef\bHK{HK}{\author{I.\ Hambleton and M.\ Kreck}
\article{Cancellation, elliptic surfaces and the topology of certain
four{--}manifolds}
\journal{J.\ Reine Angew.\ Math.}
\vol{444}
\yr{1993}
\pages 79-100.
}
\Bdef\bHKT{HKT}{\author{I.\ Hambleton, M.\ Kreck and P.\ Teichner}
\sameand{I.\ Hambleton and M.\ Kreck}{and P.\ Teichner}
\article{Nonorientable 4{--}manifolds with fundamental group 
of order 2}
\journal{Trans.\ Amer.\ Math.\ Soc.}
\yr{1994}
\vol{344}
\pages 649-665.
}
\Bdef\bHMTW{HMTW}{\author{I.\ Hambleton, R.\ J.\ Milgram,
L.\ R.\ Taylor and E.\ B.\ Williams}
\article{Surgery with finite fundamental group}
\journal{Proc.\ Lond.\ Math.\ Soc.}
\yr{1988}
\vol{56}
\pages  329-348.
}
\Bdef\bHV{H-V}{\author{J.\ C.\ Hausmann and P.\ Vogel}
\yr{1993}
\book{Geometry on Poincar\'e Spaces}
\extra{Volume 41 of {\sl Math.\ Notes.}}
\pub{Princeton Univ.\ Press}}
\Bdef\bKf{K}{\author{R.\ C.\ Kirby}
\yr{1989}
\book{Topology of 4{--}Manifolds}
\extra{Volume 1374 of {\sl Lecture Notes in Math.}}
\pub{Springer{--}Verlag}}
\Bdef\bKC{}{\author{R.\ C.\ Kirby}
\same{R.\ C.\ Kirby}
\article{Akbulut's corks and $h${--}cobordisms of smooth simply
connected 4{--}man\-i\-folds}
\journal{Turkish J.\ Math.}
\vol{20}
\yr{1996}
\pages 85-93.
}
\Bdef\bPL{PL}{\author{R.\ C.\ Kirby}
\same{R.\ C.\ Kirby}
\yr{1996}
\article{Problems in Low{--}Dimensional Topology}
\editor{W.\ H.\ Kazez (Ed.)}
\inbook{Geometric Topology, Proc.\ of the 1993 Georgia International
Topology Conference}
\pub{Amer.\ Math.\ Soc.\ \& International Press}
\pages 35-473.
\extra{AMS/IP Studies in Advanced Mathematics}
\vol{2 Part 2}
\yr{1997}
}
\Bdef\bKS{KS}{\author{R.\ C.\ Kirby and L.\ C.\ Siebenmann}
\sameand{R.\ C.\ Kirby}{and L.\ C.\ Siebenmann}
\yr{1977}
\book{Foundational Essays on Topological Manifolds, Smoothing, 
and Triangulations}
\extra{Volume 88 of {\sl Ann.\ of Math.\ Stud.}}
\pub{Princeton Univ.\ Press}}
\Bdef\bMKv{MK}{\author{M.\ Kreck}
\article{Surgery and Duality}
\preprint{Preprint, 1995}
}
\Bdef\bMKx{MKa}{\author{M.\ Kreck}
\article{$h$-Cobordisms between 4-manifolds}
\preprint{Preprint, 1995}
}
\Bdef\bMKy{MKb}{\author{M.\ Kreck}
\article{A guide to the classification of manifolds}
\preprint{This volume, 1998}
}
\Bdef\bKMa{}{\author{P.\ B.\ Kronheimer and T.\ S.\ Mrowka}
\yr{1994}
\article{Recurrence relations and asymptotics for four-manifold invariants}
\journal{Bull.\ Amer.\ Math.\ Soc}
\vol{30}
\pages 215-221.
}
\Bdef\bKMb{}{\author{P.\ B.\ Kronheimer and T.\ S.\ Mrowka}
\same{P.\ B.\ Kronheimer and T.\ S.\ Mrowka}
\yr{1995}
\article{Embedded surfaces and the structure of Donaldson's polynomial
invariants}
\journal{J.\ Differential Geom.}
\vol{41}
\pages 573-734.
}
\Bdef\bKuga{MKb}{\author{K. Kuga}
\article{Representing homology classes of $S^2 \times S^2$}
\journal{Topology}
\vol{23}
\yr{1984}
\pages 133-137.}
\Bdef\bKSh{KS}{\author{S.\ Kwasik and R.\ Schultz}
\yr{1989} 
\article{On $s${--}cobordisms of metacyclic prism manifolds}
\journal{Invent.\ Math.}
\vol{97}
\pages 523-552.
}
\Bdef\bLS{LS}{\author{R.\ K.\ Lashof and J.\ L.\ Shaneson}
\yr{1971}
\article{Smoothing four{--}manifolds}
\journal{Invent.\ Math.\ }
\vol{14}\pages 197-210.}
\Bdef\bLT{LT}{\author{R.\ K.\ Lashof and L.\ R.\ Taylor}
\yr{1984}
\article{Smoothing theory and Freedman's work on four manifolds}
\editor{I.\ Madsen and R.\ Oliver (Eds.)}
\inbook{Algebraic Topology, Aarhus, 1982}
\extra{Volume 1051 of {\sl Lecture Notes in Math.}}
\pages 271-292. 
\pub{Springer{--}Verlag}}
\Bdef\bTL{TL}{\author{T.\ Lawson}
\yr{1979}
\vol{29}
\pages 305-321.
\article{Trivializing 5{--}dimensional $h${--}cobordisms by stabilization}
\journal{Manuscripta Math.}
}
\Bdef\bLW{LW}{\author{R.\ Lee and D.\ M.\ Wilczynski}
\article{Representing homology classes by locally flat surfaces of minimal genus}
\journal{Amer. J. Math.}
\yr{1997}
\vol{119}
\pages 1119-1137.
}
\Bdef\bLL{MM}{\author{B.- H.\ Li and T.-J.\ Li}
\article{Minimal genus embeddings in $S\sp 2$-bundles over surfaces}
\journal{Math.\ Res.\ Lett.}
\vol{4}
\yr{1997}
\pages 379-394.}
\Bdef\bLLb{MM}{\author{B.- H.\ Li and T.-J.\ Li}
\article{Minimal genus smooth embeddings in $S^2 \times S^2$ and 
$\CP^2\# n\overline{\CP}^2$ with $n\leq8$}
\journal{Topology}
\vol{37}
\yr{1998}
\pages 575-594.}
\Bdef\bMM{MM}{\author{R.\ Mandelbaum and B.\ Moishezon}
\article{On the topological structure of non{--}singular algebraic
surfaces in $\bbb{\CP}^3$}
\journal{Topology}
\vol{15}
\yr{1976}
\pages 23-40.
}
\Bdef\bAstudent{}{\author{R.\ Matveyev}
\article{A decomposition of smooth simply{--}connected $h${--}cobordant
4{--}manifolds}
\journal{J.\ Differential Geom.}
\yr{1996}
\vol{44}
\pages 571-582.
}
\Bdef\bM{M}{\author{J.\ W.\ Milnor}
\yr{1964}
\article{Differential Topology}
\editor{Saaty (Ed.)}
\inbook{Lectures in Modern Mathematics, Vol. 2}
\pages 165-183. 
\pub{Wiley}}
\Bdef\bMST{MM}{\author{J. W. Morgan, Z. Szab\'o and C. H. Taubes} 
\article{A product formula for the Seiberg-Witten invariants and the
generalized Thom conjecture}
\journal{ J.\ Differential Geom.}
\vol{ 44} 
\yr{1996}
\pages 706-788.}
\Bdef\bNo{M}{\author{S.\ Novikov}
\yr{1965}
\article{Homotopy equivalent smooth manifolds I}
\journal{Translations Amer.\ Math.\ Soc.}
\vol{48}
\pages 271-396.
\pub{Wiley}}%
\Bdef\bOS{M}{\author{P. Ozsv\'ath and Z. Szab\'o}
\article{ The symplectic Thom conjecture}
\preprint{Preprint 1998}} 
\Bdef\bQa{Q1}{\author{F.\ S.\ Quinn}
\yr{1986}
\article{Isotopy of $4${--}manifolds}
\journal{J.\ Differential Geom.}
\pages 343-372.
\vol{24}
}
\Bdef\bQz{Q1}{\author{F.\ S.\ Quinn}
\yr{1997}
\article{Problems in low-dimensional topology}
\journal{Sci.\ Bull.\ Josai Univ.}
\pages 97-104.
\vol{{\rm Special issue no.} 2}
}
\Bdef\bQRt{QR}{\author{A.\ A.\ Ranicki}
\yr{1992}
\book{Algebraic $L${--}theory and Topological Manifolds}
\extra{Cambridge Tracts in Mathematics 102}
\pub{Cambridge Univ.\ Press}
}
\Bdef\bR{R}{\author{V.\ A.\ Rochlin}
\yr{1986}
\article{Quatre articles de V. A. Rohlin, IV}
\editor{L.\ Guillou and A.\ Marin (Eds.)}
\inbook{A la Recherche de la Topologie Purdue}}
\extra{Volume 62 of {\sl Prog.\ Math.}}
\pages 17-23.
\pub{Birkh\"auser}
\Bdef\bRu{xx}{\author{D.\ Ruberman}
\article{Invariant knots of free involutions of $S^4$}
\vol{18}
\journal{Topology Appl.}
\yr{1984}
\pages 217-224.
}
\Bdef\bMSa{S}{\author{M.\ G.\ Scharlemann}
\yr{1976}
\article{Transversalities theories at dimension four}
\journal{Invent.\ Math.}
\vol{33}
\pages 1-14.}
\Bdef\bMS{S}{\author{M.\ G.\ Scharlemann}
\same{M.\ G.\ Scharlemann}
\yr{1976}
\article{Constructing strange manifolds with the dodecahedral space}
\journal{Duke Math.\ J.\ }
\vol{43}
\pages 33-40.}
\Bdef\bSp{Sp}{\author{M.\ Spivak}
\yr{1967}
\article{Spaces satisfying Poincar\'e duality}
\journal{Topology}
\vol{6}
\pages 77-102.}
\Bdef\bRS{RS}{\author{R.\ Stong}
\article{A structure theorem and a splitting theorem 
for simply connected smooth 4{--}manifolds}
\journal{Math.\ Res.\ Lett.\ }
\vol 2
\yr{1995}
\pages 497-503.
}
\Bdef\bDS{Su}{\author{D.\ P.\ Sullivan}
\yr{1966}
\book{Triangulating Homotopy Equivalences}
\extra{Thesis, 1995, Princeton University}
\pub{Univ.\ Microfilms, Ann Arbor}
}
\Bdef\bTa{}{\author{C.\ H.\ Taubes}
\yr{1995}
\article{The Seiberg{--}Witten and Gromov invariants}
\journal{Math. Res. Lett.}
\vol{2}
\pages 221-238.}
\Bdef\bRT{Thom}{\author{R.\ Thom}
\yr{1954}
\article{Quelques propri\'et\'es globales 
des vari\'et\'es diff\'erentiables}
\journal{Comment.\ Math.\ Helv.\ }
\vol{28}
\pages 17-86.}
\Bdef\bCTa{1044}{%
\author{C.\ B.\ Thomas}
\yr{1970}
\article{A weak $4$-dimensional $S$-cobordism theorem}
\journal{Proc. Cambridge Philos. Soc.}
\vol{67}
\pages 549-551. 
}%
\Bdef\bCT{1044}{%
\author{C.\ B.\ Thomas}
\same{C.\ B.\ Thomas}
\yr{1995}
\article{3{--}Manifolds and $PD(3)${--}groups}
\editor{S.\ C.\ Ferry, A.\ A.\ Ranicki and J.\ M.\ Rosenberg (Eds.)}
\inbook{Novikov Conjectures, Index Theorems and Rigidity,
Volume 1 (Oberwolfach, 1993)}
\extra{Volume 226 of {\sl London Math.\ Soc.\ Lect.\ Note Ser.}}
\pages 301-308. 
\pub{Cambridge Univ.\ Press}}
\Bdef\bTh{T}{\author{W.\ P.\ Thurston}
\yr{1982}
\article{Three dimensional manifolds, Kleinian groups and 
hyperbolic geometry}
\journal{Bull.\ Amer.\ Math.\ Soc.\ }
\vol{6}\pages 357-382.}
\Bdef\bWd{W1}{\author{C.\ T.\ C.\ Wall}
\yr{1964}
\article{Diffeomorphisms of 4{--}manifolds}
\journal{J.\ London Math.\ Soc.\ }
\vol{39}
\pages 131-140.}
\Bdef\bWh{W2}{\author{C.\ T.\ C.\ Wall}
\same{C.\ T.\ C.\ Wall}
\yr{1964}
\article{On simply connected 4{--}manifolds}
\journal{J.\ London Math.\ Soc.\ }
\vol{39}
\pages141-149.}
\Bdef\bCTC{W}{\author{C.\ T.\ C.\ Wall}
\same{C.\ T.\ C.\ Wall}
\yr{1971}
\book{Surgery on Compact Manifolds}
\pub{Academic Press}}
\Bdef\bSW{}{\author{E.\ Witten}
\yr{1994}
\article{Monopoles and four-manifolds}
\journal{Math. Res. Lett.}
\vol{1}
\pages 769-796.}

\def \sectionzero{0. Review of Surgery Theory.}
\def \thmA{Theorem\ 1}
\def \thmB{Theorem\ 2}
\def \ses{3}
\def \sectionone{1. The Low Dimensional Results.}
\def \thmC{Theorem\ 4}
\def \thmtopdiff{Theorem\ 5}
\def \sectiontwo{2. Calculation of Normal Maps.}
\def \formulaB{6}
\def \thmCalc{Lemma\ 7}
\def \formulaA{8}
\def \secQR{4}
\def \sectionthree{3. Surgery Theory.}
\def \fkq{9}
\def \QRt{Proposition\ 10}
\def \tabo{Table\ 11}
\def \sectionfour{4. Computation of Stable Structure Sets.}
\def \topSfour{Theorem\ 12}
\def \topSthree{Theorem\ 13}
\def \nonE{Lemma\ 14}
\def \diffS{Theorem\ 15}
\def \diffSb{Theorem\ 16}
\def \sectionfive{5. A Construction of Novikov, Cochran \ \&\ Habegger.}
\def \afone{17}
\def \thmCH{Theorem\ 18}
\def \thmNM{Theorem\ 19}
\def \cCRA{Corollary\ 20}
\def \thmNMS{Theorem\ 21}
\def \bTHMa{Theorem\ 22}
\def \sectionsix{6. Examples.}
\def \Sfive{7}
\def \sectionseven{7. The Topological Case in General.}
\def \thmNDL{Theorem\ 23}
\def \sexp{Theorem\ 24}
\def \Ssix{8}
\def \sectioneight{8. The Smooth Case in Dimension 4.}
\def\SSS{\ifnum\truesection<100\relax{\xxx\char'170}\hskip1pt\number\truesection\fi}
\def\oddheader{\upper 4-manifolds through the eyes of surgery@ \SSS}
\def\evenheader{\upper R.\ @\upper C.\ @\upper Kirby and@
\upper L.\ @\upper R.\ @\upper Taylor@ \SSS}
\START
\topmatter
\title{A survey of 4-manifolds through the eyes of surgery.}
\authors{Robion C.\ Kirby and Laurence R.\ Taylor}
\thanks{Both authors were partially supported by the N.S.F.}
\address{Department of Mathematics\\University of California at Berkeley\\Berkeley, CA 94720\\kirby@math.berkeley.edu}
\address{Department of Mathematics\\University of Notre Dame\\Notre Dame, IN 46556\\taylor.2@nd.edu}
\dedication{To C.~T.~C.~Wall on the occasion of his sixtieth birthday.}
\endtopmatter
\font\xxx=cmsy8

Table of Contents:
{{\obeylines
\sectionzero
\sectionone
\sectiontwo
\sectionthree
\sectionfour
\sectionfive
\sectionsix
\sectionseven
\sectioneight
}}

\newsec{Review of Surgery Theory.}
\Remembersection\sectionzero
Surgery theory is a method for
constructing manifolds satisfying a given collection of homotopy
conditions.
It is usually combined with the $s${--}cobordism theorem which
constructs homeomorphisms or diffeomorphisms between two similar
looking manifolds.
Building on work of Sullivan, Wall applied these two techniques to
the problem of computing structure sets.
While this is not the only use of surgery theory, it is the aspect on which 
we will concentrate in this survey.
In dimension $4$, there are two versions, one in which one builds
topological manifolds and homeomorphisms and the second
in which one builds smooth manifolds and diffeomorphisms.
These two versions are dramatically different.
Freedman has shown that the topological case resembles the higher
dimensional theory rather closely.
Donaldson's work showed that the smooth case differs wildly from
what the high dimensional theory would predict.
Surgery theory requires calculations in homotopy theory and in low
dimensions these calculations become much more manageable.
In sections 0 and 1, we review the general theory and describe the
general results in dimensions $3$ and $4$.
In sections 2 through 6, we describe precisely what the 
high dimensional theory predicts.
Finally, we describe the current state of affairs
for the two versions in sections 7 and 8.

To begin, let $(X,\partial X)$ be a simple, 
$n${--}dimensional Poincar\'e
space whose boundary may be empty.
In particular, $X$ is homotopy equivalent to 
a finite CW complex
which satisfies Poincar\'e duality for any coefficients, 
with a twist
in the non{--}orientable case, and simple means that
there is a chain map
$$[X,\partial X]\cap\colon 
Hom_{\Z[\pi_1(X)]}(C_\ast(X),\Z[\pi_1(X)])\to
C_{n-\ast}(X)$$
which is a simple isomorphism 
between based chain complexes, \cite\bCTC.
This is the homotopy analogue of a manifold.
Let \CAT\ stand for either \TOP, 
the topological category, or
\DIFF, the differential category.
There is also the category of \PL{--}manifolds,
but it follows from the work
of Cerf, \cite\bCf, that in dimension 4 
\PL\ is equivalent to \DIFF, so we
will rarely discuss \PL\ here.
Fix a \CAT{--}manifold $L^{n-1}$ without boundary and a 
simple homotopy equivalence $h\colon L\to\partial X$.

\noindent{\bf Structure Sets:}
Define the set $\Scat(X; \rel\ h)$ as
the set of all simple homotopy equivalences of pairs,
$f\colon (M,\partial M)\to(X,\partial X)$, where $(M, \partial M)$
is a \CAT{--}manifold, and for which there exists
a \CAT{--}equivalence $g\colon L\to\partial M$ such that the
composition $L\to\partial M\to\partial X$ is homotopic to $h$;
two such, $(M_i, f_i, g_i)$ $i=0, 1$,
are deemed equal if there exists a \CAT{--}equivalence
$F\colon (M_0,\partial M_0) \to (M_1,\partial M_1)$ so that
$f_1\comp F$ is homotopic, as a map of pairs, to $f_0$, and
$F\vert_{_{\partial}}\comp g_0$ is homotopic to $g_1$.
In diagrams,
\par\noindent\hbox to\hsize{
\hskip-20pt\vtop{\hsize=3in
$$\matrix{L&\hbox to 0pt{\hss$\RA{g}$\hss}&\partial M\cr
&\hskip-10pt\hbox to 0pt{$\searrow\hskip-20pt\rzb{h}$\hss}&
\downarrow\rzb{f\vert_{\partial X}}\cr
&&\partial X\cr}$$}
\hskip-60pt\vtop{\hsize=3in
$$\matrix{M_0&\hbox to 0pt{\hss$\RA{f_0}$\hss}&X\cr
\lzb{F}\downarrow&
\hbox to 0pt{%
\hbox{$\nearrow\hskip-4pt\lower 4pt\hbox{$\rzb{f_1}$}$}\hss}
\cr
M_1\cr}$$}
\hfill}
homotopy commute.

\remark
One can use the homotopy extension theorem to tighten up the
definition: one can restrict to manifolds $M$ with $\partial M=L$
and with maps $f$ such that $f\vert_{_{\partial}}=h$;
$F\vert_{_{\partial}}$ can be required to be the identity
and the homotopy between $f_1\comp F$ and $f_0$ can be required
to be constant on $L$.
Finally, base points may be selected in each component of
$M$, $X$, $\partial M$ and $\partial X$ and all the maps and
homotopies may be assumed to preserve the base points.
This is a useful remark in identifying various fundamental groups
precisely rather than just up to inner automorphism.

The questions now are whether the set $\Scat(X; \rel\ h)$
is non{--}empty (existence) and if non{--}empty, how many elements
does it have (uniqueness).
The only 1 and 2 dimensional Poincar\'e spaces are simple
homotopy equivalent to manifolds, \cite\bEM, \cite\bEL, 
and this is conjecturally
true in dimension three, \cite\bCT.
In general, the Borel conjecture asserts that this is true for
aspherical Poincar\'e spaces in all dimensions
(see the discussion of Problem 5.29 in \cite\bPL\ and 
the articles in \cite\bFRR).

There are bundle{--}theoretic obstructions to 
$\Scat(X; \rel\ h)$ being non{--}empty.
Every Poincar\'e space has a 
stable Spivak normal fibration, \cite\bSp,
which is given by a map $\nu_X\colon X\to BG$.
This is the homotopy analogue of the stable normal bundle
for a manifold.
The space $BG$ can be thought of as 
the classifying space for stable spherical fibrations,
or as the limit of the classifying spaces of $G(m)$,
the space of homotopy automorphisms of $S^{m-1}$.
There is a map $B\CAT\to BG$ and a necessary condition for
$\Scat(X; \rel\ h)$ to be non{--}empty is that 
$\nu_X$ lift to $B\CAT$.
Given a homotopy equivalence between a \CAT{--}manifold
and a Poincar\'e  space, $X$, Sullivan, \cite\bDS,
constructs a {\sl homotopy differential},
a specific lift of $\nu_X$.
The lift to $B\CAT$ gives a stable \CAT\ bundle 
$\eta$ over $X$ and the lift gives a 
specific fibre homotopy equivalence
between the associated sphere bundle to $\eta$ 
and the Spivak normal fibration $\nu_X$.

With data as above, the Sullivan homotopy differential
gives an explicit lift of $\nu_{\partial X}$ to $B\CAT$:
a second application of this yoga gives a map
$$\N\colon \Scat(X;\rel\ h)\to \Ncat(X;\rel\ h)$$ 
where $\Ncat(X;\rel\ h)$ is the set of homotopy
classes of lifts of
$\nu_X$ to $B\CAT$ which restrict to our given 
lift over $\nu_{\partial X}$.

Boardman\and Vogt, \cite\bBV, prove that the spaces 
$B\CAT$ and $BG$ are infinite loop spaces and that
the maps $B\CAT\to BG$ are infinite loop maps.
It follows that there is a sequence of homotopy fibrations,
extending infinitely in both directions,
$$\cdots\to \CAT\to G\to \gcat \to B\CAT
\to BG\to B(\gcat)\to\cdots$$
The Spivak normal fibration is a map $\nu_X\colon X\to BG$, and
the Sullivan differential on the boundary 
gives an explicit null{--}homotopy
of $\nu_X\vert_{\partial X}$ in $B(\gcat)$ 
and so defines a map
$b\colon X/\partial X\to B(\gcat)$.

The next result follows from standard 
homotopy theory considerations:
\Ndef\thmA{Theorem}
\proclaim \thmA.
$\Ncat(X;\rel\ h)$ is non{--}empty  iff 
$b\colon X/\partial X\to B(\gcat)$ is null homotopic.
If $\Ncat(X;\rel\ h)$ is non{--}empty, the abelian group 
$\br{X/\partial X}{\gcat}$ 
acts simply{--}transitively on it.

\remark
If $X$ is already a \CAT{--}manifold, 
$\Ncat(X;\rel\ h)$ has an obvious choice of base point,
namely the normal bundle of $X$.

Given a point $x\in\Ncat(X;\rel\ h)$ and an element
$\eta\in\br{X/\partial X}{\gcat}$,
let $\eta\normact x\in\Ncat(X;\rel\ h)$ denote
the result of the action.

\CAT{--}transversality allows an 
interpretation of $\Ncat(X;\rel\ h)$
as a normal bordism theory.
We can translate this into a more geometric language
where we assume for simplicity that $\partial X=\emptyset$.
Choose a simplicial subcomplex of a high dimensional sphere, $S^N$,
which is simple homotopy equivalent to $X$.
Let $W,\partial W$
denote a regular neighborhood.
If the map $\partial W\to X$ is made into a fibration 
then the result is a spherical
fibration with fibre $S^{N-n-1}$, which is
the Spivak normal fibration; it corresponds to a 
classifying map $X\to BG$.
Note that by collapsing the complement of $W$ to a point,
we get a map from $S^N$ to the Thom space of the Spivak
normal fibration.
A lift from $BG$ to $B\CAT$ provides a fibre homotopy
equivalence from the Spivak normal fibration to the 
\CAT\ bundle over $X$, and this extends to Thom spaces.
Thus a lift from $BG$ to $B\CAT$ gives by composition
a map from $S^N$ to the Thom space of the \CAT\ bundle;
making this map transverse to the $0${--}section provides
a manifold, $M^n$, and a degree one map,
$M\to X$ covered by a \CAT\ bundle map 
from the stable normal bundle for $M$ to the given 
bundle over $X$.
Different choices change the data by a normal bordism.
Summarizing , $\Ncat(X)$
can be interpreted as bordism classes of 
degree{--}one normal maps,
that is, degree one maps $f\colon M\to X$ 
covered by a bundle map from the stable normal 
bundle of $M$ to some \CAT\ bundle over $X$.

Given a normal map $M^n\RA{h} X$, one can try to 
surger $M$ so that $h$ becomes a simple homotopy
equivalence.  
This allows one to define a surgery obstruction map 
in general,
$$\theta\colon \Ncat(X; \rel\ h)\to
\LWG sn(\Z[\pi_1(X)],w_1(X))$$
where $w_1(X)\colon \pi_1(X)\to \pm 1$ 
is the first Stiefel{--}Whitney
class of the Poincar\'e space $X$ and $\LWG sn$ is the
Wall group as defined in \cite\bCTC.
The Wall groups depend only on the group and the
first Stiefel{--}Whitney class and are 4{--}fold periodic.

In the simply connected case, the only obstruction 
in dimensions congruent to $0$ mod $4$ is the
difference in the signatures of $M$ and $X$, 
so $\LWG s0(\Z)$ is $\Z$ and the map $\theta$
is given by $(\sigma(M)-\sigma(X))/8$.
In dimensions congruent to $2$ mod $4$, do surgery
to the middle dimension,
put a quadratic enhancement on the kernel in homology 
and take the Arf invariant to get an invariant in
$\LWG s2(\Z)=\cy2$. 
The simplest example is the degree one map from 
$T^2$ to $S^2$ with stable normal map given by 
framing the stable normal bundle to $S^2$ and taking
the ``Lie framing'' of the stable normal bundle to $T^2$ 
defined as follows:
identify a normal bundle to $T^2$ 
with the product of two stable normal bundles
to $S^1$ and frame each of these with the framing 
that does not extend over $D^2$.
In odd dimensions, the obstruction is 
$0=\LWG s1(\Z)=\LWG s3(\Z)$.

If $\Scat(X;\rel\ h)\neq\emptyset$, the composite
$$\Scat(X;\rel\ h)\RA{\N}\Ncat(X;\rel\ h)\RA{\theta}
\LWG sn(\Z[\pi_1(X)],w_1(X))$$
sends every element in the structure set to 
the zero element in the Wall group.

Given $x\in \Ncat(X;\rel\ h)$, let 
$\theta_x\colon \br{X/\partial X}{\gcat}\to 
\LWG sn(\Z[\pi_1(X)],w_1(X))$
be defined by $\theta_x(\eta)=\theta(\eta\normact x)$.
Thus far, there are no dimension restrictions, but
one of Wall's fundamental results,
\ecite\bCTC{Thm 10.3 and 10.8}, is

\Ndef\thmB{Theorem}
\Ndef\ses{}
\proclaim \thmB.
If $n\geq 5$ and if
$x\in \Ncat(X;\rel\ h)$,
the following sequence is exact 
$$\Scat(X;\rel\ h)\RA{\N_x}
\br{X/\partial X}{\gcat}
\RA{\theta_x}\LWG sn(\Z[\pi_1(X)],w_1(X))\leqno(\ses)$$
in the sense that $\theta_x^{-1}(0)$ 
equals the image of $\N_x$.
If $\Scat(X;\rel\ h)\neq\emptyset$, 
there is an action of a Wall group on it:
$$\LWG s{n+1}(\Z[\pi_1(X)],w_1(X))\times 
\Scat(X;\rel\ h)\to\Scat(X;\rel\ h)$$
and $\N_x$ is injective on the orbit space.
The isotropy subgroups of this action are given 
by ``backing{--}up'' sequence (\ses), 
being careful with base point.
Specifically, if $f\colon M\to X$ is in $\Scat(X;\rel\ h)$,
let $f\times1_{[0,1]}$ be the evident map
$M\times[0,1]\to X\times[0,1]$ with 
$\partial_{f\times 1_{[0,1]}}$
being the evident homeomorphism on the boundary:
let $\N(f\times1_{[0,1]})\in\Ncat(X;\rel\ h)$
be our choice of base point, denoted $y$ below.
The isotropy subgroup of $f\colon M\to X$ is the image of 
$\theta_{y}$ in the version of (\ses)
$$\Scat(X\times [0,1]; \rel\ \partial_{f\times 1_{[0,1]}})
\RA{\N_{y}}
\br{\Sigma(X/\partial X)}{\gcat}
\RA{\theta_{y}} \LWG s{n+1}(\Z[\pi_1(X)],w_1(X))\ .$$

\newsec{The Low Dimensional Results.}
\Remembersection\sectionone
If $n<5$, sets $\barScat(X;\rel\ h)$ are defined 
below so that \thmB\ remains true if the sets 
$\barScat$ are used instead of the sets $\Scat$.
By construction there will be a map
$\psi_{\CAT}\colon \Scat(X;\rel\ h)\to\barScat(X;\rel\ h)$
and the failure of surgery in low dimensions is
the failure of $\psi_{\CAT}$ to be a bijection.

It is a fortuitous combination of calculations 
of Wall groups, the classification of manifolds
and the result that 2{--}dimensional Poincar\'e
spaces have the homotopy type of manifolds, \cite\bEM, 
\cite\bEL, that  \thmB\ holds as stated for $n=1$ and $2$.
After this remark, we restrict attention to the three
and four dimensional cases.

In dimension 3, for closed manifolds, 
it is conjectured that $\Scat(M^3)$ is a point,
\ecite\bPL{3.1$\Omega$}.
Computationally, $\barSdiff(S^3)$ is two points, $S^3$
and the Poincar\'e sphere;
however, $\barStop(S^3)$ is still one point, because
$S^3$ and the Poincar\'e sphere are 
topologically homology bordant.

In dimension 4, Freedman's work shows $\psi_{\TOP}$ is a
bijection for ``good'' fundamental groups;
Donaldson's work shows $\psi_{\DIFF}$ is not bijective 
for many 4{--}manifolds.
These points are discussed below in sections 
\Sfive\ and \Ssix.

A mantra of four{--}dimensional topology is that 
``things work after adding
$S^2\times S^2$'s'': a mantra of three{--}dimensional topology is that
``surgery works up to homology equivalence''.
The results below lend some precision to these statements.

Let us assume given $(X^3, \partial X)$ with a 
\CAT{--}homotopy structure
 $h\colon L^2\to \partial X$.
Since every 2{--}dimensional \TOP{--}manifold has a unique 
smooth structure, it is 
no loss of generality to assume $L$ is smooth.
Define $\barScat(X;\rel\ h)$ as a set of objects
modulo an equivalence relation.
Each object is a pair consisting of a \CAT{--}manifold,
$M$, and a map, $f\colon M^3\to X$, 
where $M^3$ is smooth and $f$ induces
an isomorphism in homology with coefficients 
in $\Z[\pi_1(X)]$.
Any such map has a Whitehead torsion in 
$Wh(\Z[\pi_1(X)])$ and we
further require that this torsion be $0$.
Two such objects, $M_i$, $f_i$ $i=0$, $1$, 
are deemed equivalent 
iff there exists
a normal bordism which will consist of
a \CAT{--}manifold $W^4$ with 
$\partial W= M_0\disjointunion M_1$,
a map $F\colon W\to X\times [0,1]$ extending $f_0$ and $f_1$,
a \CAT{--}bundle $\zeta$ over $X\times[0,1]$,
and a bundle map covering $F$
between the normal bundle for $W$ and $\zeta$.
In such a case, there is a well{--}defined surgery 
obstruction in $\LWG s4(\Z[\pi_1(X)],w_1(X))$ 
which we further require to be $0$.
In case \TOP{--}surgery works in dimension 4 
for $\pi_1(X)$, this condition is equivalent to 
the following more geometric statement:
if $\CAT=\TOP$, the normal bordism 
can be replaced by a topological $s${--}cobordism;
if $\CAT=\DIFF$, the normal bordism can be replaced by 
a topological $s${--}cobordism with vanishing 
stable triangulation obstruction.

We now turn to the 4{--}dimensional case.
Let us assume given $(X^4, \partial X)$ with a 
\CAT{--}homotopy structure
$h\colon L^3\to \partial X$.
Since every 3{--}dimensional \TOP{--}manifold has a unique 
smooth structure, it is 
no loss of generality to assume $L$ is smooth.
Following Wall, write $X$ as a 3{--}dimensional complex,
$\overo X\subset X$, union a single 4{--}cell.
For any integer $r> 0$, one can form the connected sum,
$X\# rS^2\times S^2$ by removing a 4{--}ball in the interior
of the top 4{--}cell.
There are maps 
$p_X\colon X\# rS^2\times S^2/\partial X \to X/\partial X$.
Define
$\wtScat{r}(X;\rel\ h)=
\{f\in\Scat(X\#rS^2\times S^2;\rel\ h)\ \vert\
\N(f)\in\image\ p_X^\ast\}$ and for uniformity, let 
$\wtScat{0}(X;\rel\ h)=\Scat(X;\rel\ h)$.
There are evident maps
$\wtScat{r}(X;\rel\ h)\to\wtScat{r+1}(X;\rel\ h)$,
so define $\barScat(X;\rel\ h)$ to be the limit.
One can define $\wtNcat{r}(X;\rel\ h)$ similarly, 
but the maps
$\wtNcat{r}(X;\rel\ h)\to \wtNcat{r+1}(X;\rel\ h)$
are isomorphisms.
We call $\barScat(X;\rel\ h)$ the {\sl stable structure set}.
\Ndef\thmC{Theorem}
\proclaim\thmC.
If $n=3$ or $4$, and if
$x\in \Ncat(X;\rel\ h)$,
the following sequence is exact 
$$\barScat(X;\rel\ h)\RA{\N_x} \br{X/\partial X}{\gcat}
\RA{\theta_x} \LWG sn(\Z[\pi_1(X)],w_1(X))\ .$$
If $\barScat(X;\rel\ h)\neq\emptyset$, 
$\LWG s{n+1}(\Z[\pi_1(X)],w_1(X))$ acts
on it and $\N_x$ is injective on the orbit space.
The isotropy subgroups are given as in \thmB.
Finally, there is a map $\psi_{\CAT}\colon 
\Scat(X;\rel\ h)\to\barScat(X;\rel\ h)$
(and, if $n=4$, $\psi^r_{\CAT}\colon\ 
\wtScat{r}(X;\rel\ h)\to\barScat(X;\rel\ h)$).

\Ddef\thmD{Addendum}
\proclaim\thmD.
If $n=4$ and if 
$f_i\colon (M_i,L)\to (X,\partial X)$, $i=0, 1$ 
are such that
$\psi(f_0)=\psi(f_1)$, there exists an s{--}cobordism, $W$,
from $M_0$ to $M_1$ which is a product over $L$, 
together with a map of pairs 
$F\colon (W,\partial W) \to (X\times[0,1],
\partial(X\times[0,1]))$
which extends $f_0$ and $f_1$ and is 
$h\times [0,1]$ on 
$L\times[0,1]\subset \partial W$.

The calculations above for the smooth and 
the topological stable structure sets can be compared
using the map $\gdiff\to\gtop$.
A second way to compare them comes from the work of 
Kirby\and Siebenmann, \cite\bKS, in high dimensions 
and proceeds as follows.
There is a function $k\colon \Stop(X;\rel\ h)
\to \br{X/\partial X}{B(\topdiff)}$
which sends $f\colon M\to X$ to the smoothing obstruction
for $M$.
The group $\br{X/\partial X}{\topdiff}$ acts on the 
smooth structure set:
an element $\eta\in \br{X/\partial X}{\topdiff}$ 
corresponds to a homeomorphism
$\hat \eta\colon M^\prime\to M$,
and let $\eta$ act on $f$ to yield 
$\eta\KSact f\colon M^\prime\RA{\hat \eta} M\RA{f}X$.
The evident relation
$\bar\eta\KSact\N(f)=\N(\eta\KSact f)$ holds,
where $\bar\eta$ denotes the composite 
$X/\partial X\RA{\eta}\topdiff\to\gdiff$.
In dimension 4, there are similar results on the 
stable structure sets thanks
to the work of Lashof\and Shaneson, \cite\bLS.
In this case $\br{X/\partial X}{B(\topdiff)}=
H^4(X,\partial X;\cy2)$
and $\br{X/\partial X}{\topdiff}=H^3(X,\partial X;\cy2)$.

\Ndef\thmtopdiff{Theorem}
\proclaim \thmtopdiff.
If $n=4$, the image of the forgetful map 
$\barSdiff(X;\rel\ h)\to\barStop(X;\rel\ h)$ is 
$k^{-1}(0)$
($k\colon\barStop(X;\rel\ h)\to H^4(X,\partial X;\cy2)$).
The group $H^3(X,\partial X;\cy2)$ 
acts on $\barSdiff(X;\rel\ h)$
and the forgetful map induces a bijection between the orbit
space and $k^{-1}(0)$.

\remark
In dimension $4$, there is another version of 
``stably \CAT\ equivalent'' that appears
sometimes in the literature.
One might say $M_1$ and $M_2$ were ``stably \CAT\ equivalent'' if
$M_1\#rS^2\times S^2$ was \CAT\ equivalent to 
$M_2\#rS^2\times S^2$.
We will rarely discuss this concept, but will say
$M_1$ and $M_2$ are {\sl weakly, stably \CAT\ equivalent}
when we do.
We say $M_1$ and $M_2$ are {\sl stably \CAT\ equivalent\/}
if there is a \CAT\ equivalence
$h\colon M_1\#rS^2\times S^2\to M_2\#rS^2\times S^2$
and a homotopy equivalence, $f\colon M_1\to M_2$,
such that $f\#r1_{S^2\times S^2}$ is homotopic to $h$.
As an indication of the difference, consider that the 
Wall group acts on our stable structure set 
(non{--}trivially in some case as we shall see below), 
whereas the top and bottom
of a normal bordism are always weakly, stably \CAT\ 
equivalent since such a bordism has a handle
decomposition with only 2 and 3 handles.
It is also easy to give examples of weakly,
stably \TOP\ equivalent,
simply connected manifolds which are not even homotopy
equivalent since there are many distinct definite forms
which become isomorphic after adding a single hyperbolic.

Kreck observes that the question of whether two manifolds
are weakly, stably \CAT\ equivalent
is a bordism question, \cite\bMKv.
More precisely, fix a map $h\colon M\to K(\pi_1(M),1)$ inducing an
isomorphism on $\pi_1$ and use the normal
bundle to get a map $h\times\nu\colon M\to K(\pi_1(M),1)\times B\CAT$.
There exists a unique class $\omega_1\in 
H^1\bigl(K(\pi_1(M),1);\cy2\bigr)$
such that $h^\ast(\omega_1)$ is the first Stiefel{--}Whitney
class of $M$.
Define $E_1(\pi_1(M),\omega_1)$ to be the homotopy fibre of the map
$K(\pi_1(M),1)\times B\CAT\RA{\omega_1\times 1+1\times w_1}
K(\cy2,1)$ and note $h\times\nu$ factors
through a map $h_1\colon M\to E_1(\pi_1(M),\omega_1)$.
The map $h_1$ induces an isomorphism on $\pi_1$: 
it induces an epimorphism on $\pi_2$ if and only if the universal
cover of $M$ is not \spin.
If the universal cover is \spin, there exists a unique class
$\omega_2\in H^2\bigl(K(\pi_1(M),1);\cy2\bigr)$ such that
$h^\ast(\omega_2)$ is the second Stiefel{--}Whitney class of $M$.
Define $E_2(\pi_1(M),\omega_1,\omega_2)$ 
as the homotopy fibre of the map
$K(\pi_1(M),1)\times B\CAT\RA{(\omega_1\times 1+1\times w_1)\times
(\omega_2\times 1+1\times w_2)}K(\cy2,1)\times K(\cy2,2)$.
Then $h$ factors through a map 
$h_2\colon M\to E_2(\pi_1(M),\omega_1,\omega_2)$
which induces an isomorphism on $\pi_1$ and an epimorphism
on $\pi_2$.
Over $E_i$, $i=1$ or $2$, there is a stable bundle coming from the map
$E_i\to B\CAT$.
One can form Thom complexes and take stable homotopy to get bordism
groups, $\Omega^{\CAT}_4(\pi_1(M),\omega_1)$ and
$\Omega^{\CAT}_4(\pi_1(M),\omega_1,\omega_2)$:
the pair $M$ and $h$ as above determine an element
$[M,h\>]\in \Omega^{\CAT}_4(\pi_1(M),\omega_1,\omega_2)$ or 
$[M,h\>]\in \Omega^{\CAT}_4(\pi_1(M),\omega_1)$ 
(depending on whether the universal cover of $M$ is \spin\ or not).
For a fixed $M$, the homotopy classes of maps $h$ 
correspond bijectively to $Out\bigl(\pi_1(M)\bigr)$, the
outer automorphism group of $\pi_1(M)$.
Define two subgroups, 
$Out\bigl(\pi_1(M),\omega_1,\omega_2\bigr)=\bigl\{
h\in Out\bigl(\pi_1(M)\bigr)\ \Big\vert\ h^\ast(\omega_1)=\omega_1\ 
\hbox{\rm and\ } \ h^\ast(\omega_2)=\omega_2\ \bigr\}$ and
$Out\bigl(\pi_1(M),\omega_1\bigr)=\bigl\{
h\in Out\bigl(\pi_1(M)\bigr)\ \Big\vert\ h^\ast(\omega_1)=\omega_1
\ \bigr\}$.

These subgroups act on the bordism groups and $M$ determines
a well{--}defined element in
$\Omega^{\CAT}_4(\pi_1(M),\omega_1,\omega_2)/
Out\bigl(\pi_1(M),\omega_1,\omega_2\bigr)$
or $\Omega^{\CAT}_4(\pi_1(M),\omega_1)/
Out\bigl(\pi_1(M),\omega_1\bigr)$
depending on whether the universal cover of $M$ is \spin\ or not.

Two manifolds $M_1$ and $M_2$ are weakly, stably \CAT\ equivalent
if and only if there exists a choice of $\omega_1$ (and $\omega_2$ if the
universal covers are \spin) such that $M_1$ and $M_2$
represent the same element in $\Omega^{\CAT}_4(\pi_1(M),\omega_1)/
Out\bigl(\pi_1(M),\omega_1\bigr)$ (or, if the universal covers are \spin, 
in $\Omega^{\CAT}_4(\pi_1(M),\omega_1,\omega_2)/
Out\bigl(\pi_1(M),\omega_1,\omega_2\bigr)$).
The proof is to construct a bordism $W^5$ between $M_1$ and $M_2$
with a map $H\colon W\to E_i$, $i=1$ or $2$ as appropriate.
Then do surgery to make $H$ as connected as possible and then
calculate that this new bordism can be built from $2$ and $3$ handles.

These bordism groups depend only on the algebraic data, 
but their calculation can be difficult.
One easy case is when $M$ is orientable ($\omega_1=0$) and
the universal cover is not \spin.
Then $\Omega^{\CAT}_4(\pi_1(M),\omega_1)$ is just the ordinary
oriented \CAT\ bordism group of $K(\pi_1(M),1)$ which is just
$H_4\bigl(K(\pi_1(M),1);\Z\bigr)\oplus \Z$ in the smooth case and
$H_4\bigl(K(\pi_1(M),1);\Z\bigr)\oplus \Z\oplus\cy2$ in the topological
case: the $\Z$ is given by the signature of $M$; the $\cy2$ is given
by the Kirby{--}Siebenmann invariant; and the element in
$H_4\bigl(K(\pi_1(M),1);\Z\bigr)$ is just $h_\ast\bigl([M\>]\bigr)$.
The action by $Out\bigl(\pi_1(M)\bigr)$ is by the identity on the
$\Z$ and the $\cy2$ and is the usual action on 
$H_4\bigl(K(\pi_1(M),1);\Z\bigr)$.

\vskip8pt
The proofs of \twolabel{\thmC}{\thmtopdiff} are relatively
straightforward given Wall's work in high dimensions.
In the 3{--}dimensional case, one simply observes 
that there are no embedding issues,
but because circles now have codimension two,
we no longer have complete control over the
fundamental group.
In the smooth case in dimension 4,
Wall, \cite\bWd, \cite\bWh, 
Cappell\and Shaneson, \cite\bCSs, and Lawson, \cite\bTL,
prove the necessary results and in the 
topological case one need
only observe that Freedman\and Quinn, \cite\bFQ,
supply the tools needed to mimic the smooth proofs.

\newsec{Calculation of Normal Maps.}
\Remembersection\sectiontwo
Given the structure of the surgery exact sequence, we
need to be able to compute the space of homotopy
classes of maps from complexes into $\gtop$ and $\gdiff$.
Standard homotopy theory tells us how to do this
in principle.

The first step in this program is to calculate the homotopy
groups of these spaces.
The surgery sequence helps in this analysis.
The L{--}groups of the trivial group 
are $\Z$, $0$, $\cy2$, $0$.

Using the ``exact sequence '' (\ses), 
the Poincar\'e conjecture and the L{--}groups show that
$\pi_i(\gtop)=\Z$, $0$, $\cy2$, $0$, 
$i\equiv 0$, $1$, $2$, $3 \hskip -4pt\pmod 4$.
Generators can be constructed as well.
In dimensions congruent to 0 mod 4,
follow Milnor, \cite\bM,  and plumb the $E_8$ form.
The boundary is a topological sphere except in dimension 4
where it is the Poincar\'e homology sphere.
Cone the boundary or use Freedman, \cite\bFa, to complete
to a closed manifold, denoted $E_8$,
and construct a normal degree one map to the sphere. 
In dimensions congruent to 2 mod 4, 
follow a similar process.
Plumb two tangent bundles to
$S^{2k+1}$.  The boundary is a homotopy sphere.
Cone the boundary to get a PL manifold, $M^{4k+2}$,
and a degree one
map $f\colon M\to S^{4k+2}$.
This map can be made into a normal map 
so as to have non{--}zero 
surgery obstruction
(already done in dimension 2 above as a map $T^2\to S^2$).
See e.g.\ Browder, \cite\bBB, \S V.

One can do a similar analysis on $\pi_i(\gdiff)$ 
except now the Poincar\'e
conjecture fails in high dimension.
Still, $\pi_i(\gdiff)=\pi_i(\gtop)$ for $i<8$, 
although the map
$\pi_4(\gdiff)\to\pi_4(\gtop)$ is 
multiplication by 2 
(Rochlin's theorem, \cite\bR, or \cite\bKf).
Purists will quibble that the results used above
require the calculations they are quoted to justify, 
but the quoted results are correct and proved 
ten years before Freedman's work by Sullivan, \cite\bDS, 
Kirby\and Siebenmann, \cite\bKS.

The first two stages of a Postnikov decomposition
for $\gcat$ are 
$$K(\Z,4) \to \gcat\to K(\cy2,2)\ .$$
Rochlin's theorem shows that normal maps over $S^4$
have surgery obstruction divisible by 16;
on the other hand, there is a normal
map $M=\CP^2\#8\>\overline{\CP}^{\,2}\to \CP^2$,
defined as follows.
The cohomology class $(3,1,\cdots,1)$ determines a 
degree one map, $f\colon M\to \CP^{\,2}$.
Note $7$ times the Hopf bundle pulls back via
$f$ to the normal bundle of $M$.
As Sullivan observes, this means the first 
$k${--}invariant of $\gdiff$  is non{--}zero.
This $k${--}invariant lives in $H^5(K(\cy2,2);\Z)=\cy4$,
\cite\bBT;
$\gdiff$ is an $H${--}space so its $k$ invariants are 
primitive
\footnote{$^1$}{A primitive in the cohomology of 
an $H${--}space,
$m\colon Y\times Y\to Y$, is a cohomology class $y$ 
such that $m^\ast(y)=1\times y+y\times 1$.}.
In $H^5(K(\cy2,2);\Z)$ only $0$ and $2$ are primitives,
\cite\bBT.
Hence the first $k${--}invariant for $\gdiff$ is 
$2$, which as a cohomology operation is $\delta Sq^2$,
the integral Bockstein of the second Steenrod square.
Freedman's construction of the $E_8$ manifold shows that
the first $k${--}invariant of $\gtop$ is trivial.
(Again, Kirby\and Siebenmann had already shown this result, 
but the above makes a nice justification for the result.)

The next $k$ invariant for both $\gdiff$ and
$\gtop$ is trivial, so in particular there are maps
\Ndef\formulaB{}
$$\eqalign{
\gtop\to &\ K(\cy2,2)\times
\hphantom{\scriptscriptstyle\delta Sq^2}
K(\Z,4)\cr
\gdiff \to &\ K(\cy2,2)\times_{_{\delta Sq^2}}
K(\Z,4)\cr
}\leqno(\formulaB)$$
which are 5{--}connected.

The first $k${--}invariant of $\Omega\gdiff$ 
is the composition $\Omega(\delta)\comp\Omega(Sq^2)$ and
$\Omega(Sq^2)=0$.
This remark is useful in computing
$\br{\Sigma Y}{\gdiff}=\br{Y}{\Omega(\gdiff)}$.

Having computed the first $k${--}invariants for these
spaces, we want to extract explicit calculations of the
groups $\br Y{\gcat}$ for $Y$ a 4{--}complex as well as
a calculation of the map induced by the map
$\gdiff\to\gtop$.
There is a class $k\in H^4(B\TOP;\cy2)$, 
the stable triangulation
obstruction, which restricts to a class, 
$k\in H^4(\gtop;\cy2)$.
This class certainly vanishes when restricted 
to $\gdiff$ and
we wish to identify it in $H^4(\gtop;\cy2)$.
Let $f\colon M^4\to N^4$ be a normal map.
By \thmA, $f$ corresponds to a map $\hat f\colon N\to\gtop$
and the composite $N\RA{\hat f}\gtop\to B\TOP$
determines a bundle $\zeta$ over $N$ such that
$\nu_N\oplus\zeta$ pulls back via $f^\ast$ to $\nu_M$.
Then $k(\nu_M)=k(\zeta)+k(\nu_N)$ so $\hat f^\ast(k)$
is the difference of the triangulation obstructions for
$M$ and $N$.
Now $H^4(\gtop;\cy2)=\cy2\oplus\cy2$ generated by
$\iota_2^2$ and
$(\iota_4)_2$.
Here $\iota_2\in H^2(K(\cy2,2);\cy2)$ and
$\iota_4\in H^4(K(\Z,4);\Z)$ are generators and
$(\iota_4)_2$ denotes the generator of $H^4(K(\Z,4);\cy2)$.
By examining the normal maps,
$\chern\to \CP^2$ (where $\chern$ is 
Freedman's Chern manifold, \cite\bFa) and $E_8\to S^4$
one sees $$k=\iota_2^2+(\iota_4)_2\ .$$
One can further see that if $\hat f(k)=0$, 
then the map $N\to\gtop$
factors through a map $N\to \gdiff$.

Let $X$ be a connected 4{--}dimensional Poincar\'e space.
The maps in (\formulaB) induce natural equivalences of
abelian groups,
$$\eqalign{
\br{X/\partial X}{\gtop} &= H^2(X,\partial X;\cy2)
\oplus H^4(X,\partial X;\Z)\cr
\br{\Sigma(X/\partial X)}{\gtop} &= H^1(X,\partial X;\cy2)
\oplus H^3(X,\partial X;\Z)\cr
}$$
The calculations for $\gdiff$ look similar:
$$\eqalign{
0\to H^4(X,\partial X;\Z)\to\br{X/\partial X}{\gdiff}
& \to H^2(X,\partial X;\cy2)\to0\cr
\br{\Sigma(X/\partial X)}{\gdiff}
&= H^1(X,\partial X;\cy2)\oplus H^3(X,\partial X;\Z)\ .\cr
}$$
In general, the exact sequence for $\gdiff$ is not split.
To describe the result, let ${\cal H}^2(X,\partial X)$
denote the kernel
of the homomorphism given by the cup square, 
$H^2(X,\partial X;\cy2) \to H^4(X,\partial X;\cy2)=\cy2 $.
Note ${\cal H}^2(X,\partial X)=H^2(X,\partial X;\cy2)$ iff
$v_2(X)=0$ where $v_2$ denotes the second Wu class 
of the tangent bundle.
\Ndef\thmCalc{Lemma}
\Ndef\formulaA{}
\proclaim \thmCalc.
For $X$ a connected 4{--}dimensional Poincar\'e space with boundary,
$$\eqalign{
\hbox to 40pt{$(\ast)\hfil$}\br{X/\partial X}{\gdiff} 
=& H^2(X,\partial X;\cy2)\oplus H^4(X,\partial X;\Z)
\hskip26pt{\sl if\ } v_2(X)=0\cr
\hbox to 40pt{$(\ast\ast)\hfil$}\br{X/\partial X}
{\gdiff} =& 
{\cal H}^2(X,\partial X)\oplus \cases {\Z &if $w_1(X)=0$
and $v_2(X)\neq0$\cr
\cy4&if $w_1(X)\neq0$ and $v_2(X)\neq0$\cr}\cr
}\leqno(\formulaA).$$
The splitting in case ($\ast\ast$) depends on the 
choice of an element
$x\in H^2(X,\partial X;\cy2)$ of odd square.
The map of $\br{X/\partial X}{\gdiff}$ into 
$\br{X/\partial X}{\gtop}=H^2(X,\partial X;\cy2)
\oplus H^4(X,\partial X;\Z)$
in case ($\ast$) is just an isomorphism on $H^2$ 
and multiplication by 2
on $H^4$ and in case ($\ast\ast$) it is inclusion 
on ${\cal H}^2$ and sends
the generator of the $\Z$ (respectively $\cy4$) 
to $(x,1)$ where $1$ denotes
a generator of $H^4(X,\partial X;\Z)=\Z$ 
(respectively $\cy2$).

\remark
For 3{--}dimensional Poincar\'e spaces, the map
$\gcat\to K(\cy2,2)$ induces an isomorphism,
$[X/\partial X,\gcat]\to H^2(X,\partial X;\cy2)$.

A proof of \thmCalc\ can be constructed along the following
lines.
A diagram chase shows that 
$\br{X/\partial X}{\gdiff} \to \br{X/\partial X}{\gtop}$
is injective whenever $X$ is orientable: the image
is the kernel of $k$.
Another diagram chase shows that every element in 
${\cal H}^2(X,\partial X)\subset H^2(X,\partial X;\cy2)$
lifts to an element of order 2 in
$\br{X/\partial X}{\gdiff}$
and any lift of an
element of odd square to $\br{X/\partial X}{\gdiff}$ 
has infinite order.
This is formula \formulaA\ in the orientable case.

Assume $X$ in non{--}orientable.
If $\partial X\neq\emptyset$, let $D(X)$ denote 
the double of $X$.
Since $X\subset D(X)\to X/\partial X$ is a cofibration
and since
the inclusion $X\subset D(X)$ is split,
the case with  boundary follows from the closed case.
From Thom, \cite{\bRT}, there exists a smooth manifold
and a map $f\colon M^4\to X$ which is an isomorphism on 
$H_4(\>\vrule width 6pt depth .1pt height 0pt \>;\cy2)$.
It then follows that $f^\ast$ is an isomorphism on 
$H^4(\>\vrule width 6pt depth .1pt height 0pt \>;\Z)$ 
and an injection on
$H^2(\>\vrule width 6pt depth .1pt height 0pt \>;\cy2)$.
Hence $f^\ast$ is injective on 
$\br{\>\vrule width 6pt depth .1pt height 0pt \>}{\gdiff}$
so we may assume $X$ is a smooth manifold.
Every 2{--}dimensional homology class is represented by an
embedded submanifold, $F\subset X$,
and hence the Poincar\'e dual is the
pull back of a map $X\to T(\eta)$, 
where $\eta$ is a 2{--}plane bundle over $F$.
A diagram chase reduces the proof of \thmCalc\ 
to the calculation for $T(\eta)$.
Smashing the part of $F$ outside a disk to a point 
gives a map $F\to S^2$,
and there is a bundle $\nu$ over $S^2$ with a map
$T(\eta)\to T(\nu)$.
The bundle $\nu$ is classified by an integer,
its Euler class,
and it follows from the oriented result above that
$$\br{T(\nu)}{\gdiff}=\cases{\Z& if $\chi(\nu)$ is odd\cr
\Z\oplus\cy2 & if $\chi(\nu)$ is even\cr}\ ,$$
where the $\cy2$ in case $\chi(\nu)$ odd maps onto 
$H^2(T(\nu);\cy2)=\cy2$.
This implies \thmCalc\ in general.

The remaining question concerning normal maps 
is whether $\Ncat(X;\rel\ h)$ is empty or not:
homotopy theory says that the Spivak normal bundle
plus the lift over $\partial X$ defines a map
$X/\partial X\to B(\gcat)$.
In the \TOP\ case, $\br{X/\partial X}{B(\gtop)}=
H^3(X,\partial X;\cy2)$.
The class $g_3\colon BG\to B(\gtop)\to K(\cy2,3)$
was defined by Gitler\and Stasheff, \cite{\bGS}.
One can show that $g_3$ evaluates non{--}trivially on
$\pi_3(BG)=\cy2$.
The generator of $\pi_3(BG)$ 
corresponds to the generator of the stable 2{--}stem, 
since $\pi_{k+1}(BG)$ is isomorphic to the stable 
$k${--}stem for all $k$.
This in turn can be understood via the Pontrjagin{--}Thom
construction as a map from $S^4$ to $S^2$ 
with the inverse image of a point
being $T^2$ with the ``Lie group framing''.

Hambleton\and Milgram, \cite{\bHM}, construct a 
non{--}orientable Poincar\'e space with $g_3\neq 0$.
Using the Levitt{--}Jones{--}Quinn Poincar\'e 
bordism sequence,
\ecite{\bHV}{4.5 p.90}, one can analyze
this situation in the oriented case as well. 
One sees that $g_3$ always vanishes in the closed, 
orientable 4{--}dimensional case, as well as in
the 3{--}dimensional case.

\newsec{Surgery Theory.}
\remembersection\secQR
\Remembersection\sectionthree
The Quinn{--}Ranicki theory, \cite\bQRt,
of the assembly map can be used to decouple
the surgery theory from the specifics of the Poincar\'e space $X$.
More precisely, this section defines groups which
depend only on the fundamental group, the orientation,
the fundamental groups of the boundary
and the image of the fundamental class of $X$
in the homology of the fundamental group rel the 
fundamental group(s) of the boundary.
One of these groups will be a quotient of $L_5$ and
will act freely on the structure set so that the
quotient injects into the set of normal maps.
Another acts freely on the smooth structure set so that the
orbit space injects into the topological structure set.
Yet another gives a piece of the set of normal maps.
The results of Quinn and Ranicki
are one of the major developments
in general surgery theory
and provide the following description of the surgery obstruction map.

A Poincar\'e space with a lift of its Spivak normal
fibration to $BTOP$ acquires a fundamental class
in a twisted, $n${--}dimensional 
extraordinary homology theory, ${\bf L}^0$.
The theory ${\bf L}^0$ is a ring theory and there is a theory,
 $\lct$, so that 
$\br{X/\partial X}{\gtop}$ is the $0${--}th cohomology
group for $\lct${--}theory and $\cap D$ is just the usual
Poincar\'e duality isomorphism 
given by cap product with the fundamental class,
$\cap [X]\colon \br{X/\partial X}{\gtop}\to \lct_n^{w_1(X)}(X)$.
The map classifying the universal cover, $u\colon X\to B\pi$
induces a map 
$u_\ast\colon  \lct_n^{w_1(X)}(X)\to  \lct_n^{w_1}(B\pi)$.
There is a map $A$, the assembly map,
$$A_{\pi_1,w_1}\colon
\lct_n^{w_1}(B\pi_1)\to
\LWG sn(\zpi 2{\pi_1},w_1)\ .$$

The composite 
$\alpha=A_{\pi_1(X),w_1(X)}\comp u_\ast\comp (\cap [X])$,
$$
\br{X/\partial X}{\gtop}
\RA{\cap [X]} \lct_n^{w_1(X)}(X)\RA{u_\ast}
\lct_n^{w_1}(B\pi)\RA{A}
\LWG sn(\zpi 2{\pi_1(X)},w_1(X))$$
is related to surgery via the following formula:
let $x\in\Ncat(X;\rel\ h)$ be a chosen basepoint; 
then for any $\eta\in\br{X/\partial X}{\gtop}$
$$\alpha(\eta)=\theta(\eta\normact x)-\theta(x)\ .$$
If $X$ has the homotopy type of a manifold, $x$ can 
be chosen so that $\theta(x)=0$ and in general this 
approach divides the problem into a homotopy part and 
an algebraic part, $A_{\pi,w_1}$.
Since $A_{\pi,w_1}$ is a purely algebraic object, 
one can attack its analysis via algebra or via topology 
by using known structure set calculations.
As an example, the Poincar\'e conjecture for $n\geq5$ says 
$\Stop(S^n)$ has one point and 
one sees that the assembly map for the trivial group 
must be an isomorphism for this to work.

For analyzing the 4{--}dimensional case, 
we need to understand
$\lct_4$ and $\lct_5$; the 3{--}dimensional case 
requires that we also understand $\lct_3$.
The Atiyah{--}Hirzebruch spectral sequence for 
$\lct_\ast$ collapses for $\ast<8$ since all the 
differentials are odd torsion:
hence, for any space $Y$ and $w_1\in H^1(Y;\cy2)$,
$$\matrix{
\lct_\ast^{w_1}(Y)= 0,\  \ast\leq 1\hfill
&\lct_3^{w_1}(Y)=\hbox to 1.6in{\hfill$H_1(Y;\cy2)$}\hfill\cr
\lct_2^{w_1}(Y)=H_0(Y;\cy2)\hfill
&\lct_4^{w_1}(Y)=\hbox to 1.6in{\hfill
$H_0(Y;\Z^{w_1})\oplus H_2(Y;\cy2)$}
\hfill\cr
&\lct_5^{w_1}(Y)=\hbox to 1.6in{\hfill
$H_1(Y;\Z^{w_1})\oplus H_3(Y;\cy2)$}
\hfill\cr
}$$

\Ndef\fkq{}
Define ${\cal K}_n(\pi,w_1)$ and $Q_n(\pi,w_1)$ so as to make
$$0\to {\cal K}_n(\pi,w_1)\to
\lct_n^{w_1}(B\pi)\RA{A_{\pi,w_1}}
 \LWG sn(\zpi 2{\pi},w_1)\to Q_n(\pi,w_1)\to0\leqno(\fkq)$$
exact.

The sequences (\fkq) for various $n$ clearly only depend
on $\pi$ and $w_1$.
The groups needed for calculating the stable structure sets,
$\barScat(X;\rel\ h)$ should have the
$\lct_n^{w_1}(B\pi)$ replace by
$\lct_n^{w_1}(X)$ using the map $u_\ast$.
In the 3{--}dimensional case, $u_\ast$ is an isomorphism;
for the 4{--}dimensional case
$u_\ast$ is still an epimorphism.
For the dimensions considered here, the 5{--}dimensional
case is only needed to compute the action of the
$L${--}group on the structure set.
We want to identify the quotient group of $\LWG s5$
which acts freely, but $Q_5$ is usually too small.
The map $H_1(X;\Z^{w_1})\to H_1(B\pi;\Z^{w_1})$ is an
isomorphism, but the map $H_3(X;\cy2)\to H_3(B\pi;\cy2)$
needs to be analyzed.
The boundary of $X$ may have several components,
each with its own fundamental group:
let $\cup B\pi_1(\partial X)$ be {\sl notation\/}
for the disjoint union of the classifying spaces for
the fundamental groups of the various components
of the boundary.
There is a class 
$$D_X\in H_4(B\pi_1(X),\cup B\pi_1(\partial X);\cy2)$$
which is the image of the fundamental class of the Poincar\'e space.
Cap product with $D_X$ defines a homomorphism,
$\cap D_X\colon H^1(B\pi_1(X),\cup B\pi_1(\partial X);\cy2)
\to H_3(B\pi_1(X);\cy2)$ which is the image of $u_\ast$.
Let $\bar\lct_5^{w_1}(B\pi)=H_1(B\pi;\Z^{w_1})\oplus
H^1(B\pi_1(X),\cup B\pi_1(\partial X);\cy2)$ and let
$\bar\lct_5^{w_1}(B\pi)\to\lct_5^{w_1}(B\pi)$ be the map
which is the identity on $H_1$ and $\cap D_X$ on $H^1$.
Let $\bar A_{\pi,w_1}\colon
\bar\lct_5^{w_1}(B\pi)\to\lct_5^{w_1}(B\pi)
\RA{A_{\pi,w_1}} \LWG s5(\zpi1\pi,w_1)$, and define
$\bar{\cal K}_5(\pi,w_1,D_X)$ and $\bar Q_5(\pi,w_1,D_X)$
as the kernel and cokernel of $\bar A_{\pi.\,w_1}$.
Define
$$\bar\gamma(\pi_1(X),w_1(X),D_X)=H_1(B\pi_1(X);\cy2)/
p_1(\bar{\cal K}_5(\pi_1(X),w_1(X),D_X))$$
where 
$$p_1\colon H_1(B\pi_1(X);\Z^{w_1(X)})\oplus
H^1(B\pi_1(X),\cup B\pi_1(\partial X);\cy2)
\to H_1(B\pi_1(X);\cy2)$$ 
denotes the evident projection.

As we shall see, this $\bar\gamma$ describes the
difference between the \TOP\ and \DIFF{--}structure sets.
Define two pairs of groups depending only on
$\pi$ and $w_1$ so that
$$0\to\hat{\cal K}_5(\pi,w_1)\to
H_1(B\pi;\Z^{w_1})\to \LWG s5(\Z\pi,w_1)\to
\hat Q_5(\pi,w_1)\to0$$
is exact and define
$\hat\gamma(\pi,w_1)=H_1(B\pi_1(X);\cy2)/
p_1(\hat{\cal K}_5(\pi_1(X),w_1(X)))$
and
$\gamma(\pi,w_1)=H_1(B\pi_1(X);\cy2)/
p_1({\cal K}_5(\pi_1(X),w_1(X)))$.

\Ndef\QRt{Proposition}
\proclaim \QRt.
There are epimorphisms
$\hat\gamma\to\bar\gamma\to\gamma$ and
$\hat Q_5\to \bar Q_5\to Q_5$.
\item {1.}If $\LWG s1(\Z\pi,w_1)=0$, then
$\hat Q_5=\bar Q_5= Q_5=0$ and 
\item{} $\displaystyle\hat\gamma=\bar\gamma=\gamma=
\cases{0&if $w_1$ is trivial\cr \cy2& otherwise\cr}$
\item{2.} If $H_3(B\pi;\cy2)=0$, or if $D_X=0$,  or if 
$H^1(B\pi_1(X),\cup B\pi_1(\partial X);\cy2)=0$,
or if $\LWG s1(\Z\pi,w_1)$ has no 2{--}torsion,
then $\hat Q_5\to \bar Q_5$ and $\hat\gamma\to\bar\gamma$ 
are isomorphisms.

Two of the big conjectures in surgery theory have 
direct implications here.
The Novikov conjecture says that the $A_{\pi,w_1}$
are injective after tensoring with $\Q$.
The Borel conjecture implies that, if $B\pi$ 
is a finite Poincar\'e complex, 
then $A_{\pi,w_1}$ is split injective.
Both of these conjectures are known to be true in 
many examples.

\Ndef\tabo{Table}
Here is a table of some sample calculations.
In all cases of \tabo, \QRt\ applies:
moreover, the Whitehead group vanishes and
${\cal K}_2=0$ for all the listed groups:
$Q_2=0$ for all the listed groups except $\Z\oplus\Z$.
The displayed calculations are drawn from many sources.\null
\vskip4pt

{{
\def\CY#1{\Z/#1}
\parskip=0pt
\parindent=0pt

\baselineskip=8pt
\HSIZE=4.95in
\global\advance\HSIZE by4pt
\edef\pp{\par\hbox to 0pt{%
\hskip-4pt\advance\HSIZE by 4pt\vrule width \HSIZE depth -4pt height 4.1pt\hss}\par}
\edef\PP{\par\hbox to 0pt{%
\hskip-4pt\advance\HSIZE by 4pt\vrule width \HSIZE depth -4pt height 5pt\hss}\par}
\edef\pP{\par\hbox to 0pt{%
\hfil}\par}

\HAL=8pt
\setbox0=\hbox{$\CY2$}
\def\lz#1{\hskip\wd0\hbox to0pt{\hss$#1$}}
\def\rz#1{\hbox to0pt{$#1$\hss}\hskip\wd0}

\setbox0=\vtop{\hsize=6in%
\halign{$#$\hfil&\hskip\HAL\hfil$#$\hfil&
\hskip\HAL\hfil$#$\hfil&
\hskip\HAL\hfil$#$\hfil&
\hskip\HAL\hfil$#$\hfil&
\hskip\HAL\hfil$#$\hfil&
\hskip\HAL\hfil$#$\hfil&
\hskip\HAL\hfil$#$\hfil\cr\noalign{\pp}
\pi&\{e\}&\CY2&\CY2&\Z&\Z&\lz{\CY2}\oplus\rz\Z&\Z\oplus\Z
\cr\noalign{\pp}
w_1&0&iso.&0&0&epi.&0&0
\cr\noalign{\PP}
L_0(\Z\pi,w_1)&\Z&\CY2&\Z\oplus\Z\hphantom{/2}&\Z&\CY2&
\Z\oplus\Z\oplus\CY2&\lz\Z\oplus\rz{\CY2}
\cr\noalign{\pp}
L_1(\Z\pi,w_1)&0&0&0&\Z&0&\Z\oplus\Z&\Z\oplus\Z
\cr\noalign{\pp}
L_2(\Z\pi,w_1)&\CY2&\CY2&\CY2&\CY2&\CY2&\CY2&\lz\Z\oplus\rz{\CY2}
\cr\noalign{\pp}
L_3(\Z\pi,w_1)&0&0&\CY2&\CY2&\CY2&
\CY2\oplus\CY2&\CY2\oplus\CY2
\cr\noalign{\PP}
\lct_4^{w_1}(B\pi)&\Z&\CY2\oplus\CY2&
\Z\oplus{\CY2}&\Z&\CY2&\Z\oplus\CY2\oplus\CY2&
\lz\Z\oplus\rz{\CY2}
\cr\noalign{\pp}
H_1(B\pi;\Z^{w_1})&0&0&\CY2&\Z&0&\lz\Z\oplus\rz{\CY2}&\Z\oplus\Z
\cr\noalign{\pp}
\lct_3^{w_1}(B\pi)&0&\CY2&\CY2&\CY2&\CY2&
\CY2\oplus\CY2&\CY2\oplus\CY2
\cr\noalign{\PP}
{\cal K}_4&0&\CY2&\CY2&0&0&\CY2&0
\cr\noalign{\pp}
Q_4&0&0&\Z&0&0&\Z&0
\cr\noalign{\pp}
{\cal K}_3&0&\CY2&0&0&0&0&0
\cr\noalign{\pp}
Q_3&0&0&0&0&0&0&0
\cr\noalign{\pp}
\hat\gamma&0&\CY2&0&\CY2&\CY2&\CY2&\CY2\oplus\CY2
\cr\noalign{\pp}
\hat Q_5&0&0&0&0&0&\Z&0
\cr\noalign{\pp}
}}
\global\dEP=\dp0\relax
\hbox{\hskip.3in%
\vbox{\hsize=5.2in\box0
\vskip -5pt\vbox to 0pt{\vss%
\hskip-4pt\vrule width .1pt height \dEP depth 0pt%
\hskip .9in\vrule width .1pt height \dEP depth 0pt%
\hskip .35in\vrule width .1pt height \dEP depth 0pt%
\hskip .75in\vrule width .1pt height \dEP depth 0pt%
\hskip .6in\vrule width .1pt height \dEP depth 0pt%
\hskip .35in\vrule width .1pt height \dEP depth 0pt%
\hskip .35in\vrule width .1pt height \dEP depth 0pt%
\hskip 1in\vrule width .1pt height \dEP depth 0pt%
\hskip .75in\vrule width .1pt height \dEP depth 0pt%
}
}}
}}
\vskip 10pt
\centerline{\bf \tabo}
\vskip10pt
There are some results of a general nature which follow
from naturality and the above calculations.
If $w_1$ is trivial, ${\cal K}_4$ is a subgroup of $H_2(B\pi;\cy2)$ and
${\cal K}_2={\cal K}_3=0$.
If $w_1$ is non{--}trivial, ${\cal K}_4$ is at most 
$H_2(B\pi;\cy2)\oplus\cy2$
and ${\cal K}_2=0$.
More calculations for finite groups can be deduced from
\cite{\bHMTW}.

\newsec{Computation of Stable Structure Sets.}
\Remembersection\sectionfour

The stable \TOP{--}structure sets can now be ``computed'' .
First of all there is nothing to do if 
$\Ntop(X;\rel\ h)=\emptyset$ 
so assume it is non{--}empty 
(as it always is in the 3{--}dimensional and the
orientable 4{--}dimensional cases) and let
$$\hat\theta\colon \Ntop(X;\rel\ h)\RA{\theta}
\LWG sn(\zpi2{\pi_1(X)},w_1(X))
\to Q_n(\pi_1(X),w_1(X))\ .$$
By the surgery theory in the last section,
the image of $\hat\theta$ 
is a single point, denoted $\hat\theta(X,\rel\ h)$.

\Ndef\topSfour{Theorem}
\proclaim \topSfour: \TOP{--}structures for $\bf n=4$.
$\barStop(X;\rel\ h)\neq\emptyset$ 
iff $\hat\theta(X,\rel\ h)$ is the $0$ element in $Q_4$.
If the stable structure set is non{--}empty,
$\bar Q_5(\pi_1(X),w_1(X),D_X)$ acts freely on it.
Choose a base point $\ast$ in it.
Then $\N_{\N({\ast})}\colon \barStop(X;\rel\ h)\to
[X/\partial X;\gtop]$ induces a
bijection between the orbit space and the subgroup
of $H_2(X;\cy2)$ which maps onto 
${\cal K}_4(\pi_1(X), w_1(X))\subset H_2(B\pi_1(X);\cy2)$.

\remark
If $\pi_1(X)$ is trivial, $N_{N(\ast)}$ identifies
$\barStop(X;\rel\ h)$ with $H_2(X;\cy2)$.
In \cCRA\ below, it is shown that although the structure set can be
large there is always just one or two distinct manifolds
in it.

The 3{--}dimensional case is even easier.

\Ndef\topSthree{Theorem}
\proclaim \topSthree: \TOP{--}structures for $\bf n=3$.
$\barStop(X;\rel\ h)\neq\emptyset$ 
iff $\hat\theta(X,\rel\ h)$ is the $0$ element in $Q_3$.
If the stable structure set is non{--}empty,
$Q_4(\pi_1(X),w_1(X))$ acts freely on it.
Choose a base point $\ast$ in it.
Then $\N_{\N({\ast})}$ induces a
bijection between the orbit space and 
${\cal K}_3(\pi_1(X), w_1(X))$.

To analyze the stable smooth structure set, 
we need good criteria to see if it is non{--}empty.
Assuming $\barStop(X;\rel\ h)\neq\emptyset$,
the stable smoothing obstruction is a function
$k\colon \barStop(X;\rel\ h)\to H^4(X,\partial X;\cy2)$ 
and $\barSdiff(X;\rel\ h)\neq\emptyset$
iff $k^{-1}(0)\neq\emptyset$ (see \thmtopdiff).
In particular, it is non{--}empty in the
3{--}dimensional case.
In the simply connected, 4{--}dimensional, case, 
Freedman, \cite\bFa, argues that $k$ is constant iff
$X$ is \spin, and he constructs examples 
where the constant is 0 and others
where the constant is 1.
In the non{--}simply connected case, $v_2(X)=0$
still implies $k$ constant, but life is more complicated
when $v_2(X)\neq0$.
To describe the situation, let $\widetilde X\to X$
denote the universal cover.
If $\widetilde X$ is not \spin, then $k$ is not constant.
If $\widetilde X$ is \spin, then there exists a
unique class $v\in H^2(B\pi_1(X);\cy2)$ such that
$u^\ast(v)=v_2(X)$ under the map
$u\colon X\to B\pi_1(X)$ which classifies 
the universal cover.
Evaluation yields a map
$\cap v\colon H_2(B\pi_1(X);\cy2)\to\cy2$.

\Ndef\nonE{Lemma}
\proclaim\nonE.
$k$ is constant iff $\widetilde X$ is \spin\ and
${\cal K}_4(\pi_1(X),w_1(X))\subset \ker(\cap v)$.

\remark
If $\pi$ is finitely presented,
any classes $w\in H^1(B\pi,\cy2)$ and 
$v\in H^2(B\pi;\cy2)$ 
can be $w_1$ and $v_2$ for a
manifold with universal cover \spin.
Hence, as soon as
${\cal K}_4(\pi,w_1)\neq H_2(B\pi;\cy2)$, there are 
examples of manifolds with constant $k$ 
for which $v_2\neq0$.
From \tabo, $\Z\oplus\Z$ is such a group.
For an explicit example, recall 
$\CP^2\#\bar{\CP}^2\to S^2$
is a 2{--}sphere bundle with $w_2\neq 0$.  
Pull this bundle back over the degree one map 
$T^2\to S^2$ and
let $M^4$ denote the total space.
Then $\widetilde M$ is \spin, but $M$ is not: 
nevertheless, $k$ is constant.

\Ndef\diffS{Theorem}
\proclaim \diffS: \DIFF{--}structures for ${\bf n=4}$.
If $k^{-1}(0)\subset\barStop(X;\rel\ h)$ is non{--}empty,
$\barSdiff(X;\rel\ h)\neq\emptyset$.
The group 
$\bar\gamma(\pi_1(X),w_1(X),D_X)$ 
acts freely on $\barSdiff(X;\rel\ h)$;
the orbit space is the subset $k^{-1}(0)$.

\Ndef\diffSb{Theorem}
\proclaim \diffSb: \DIFF{--}structures for ${\bf n=3}$.
If $\barStop(X;\rel\ h)$ is non{--}empty, then
$\rho\colon\barSdiff(X;\rel\ h)\to\barStop(X;\rel\ h)$ is onto.
If $w_1(X)^2=0$, $\rho$ is 2 to 1;
if $w_1(X)^2\neq0$, $\rho$ is a bijection.

\remark
By Poincar\'e duality
$H_1(B\pi_1(X);\cy2)=H_1(X;\cy2)=H^3(X,\partial X;\cy2)$
so the action of $\bar\gamma$ gives an action of 
$H^3(X,\partial X;\cy2)$ on $\barSdiff(X;\rel\ h)$
which is the Kirby\and Siebenmann action as extended
by Lashof\and Shaneson to dimension 4.
In dimension 3, $\cy2$ acts on $\barSdiff(X;\rel\ h)$
by forming the connected sum with the Poincar\'e sphere.
If $w_1(X)^2\neq0$, this action is trivial, otherwise it is free.

The proofs of these results are fairly straightforward.
The \TOP{--}results follow from the sequence (\ses) 
for \TOP\ and the results from \S\secQR.
The \DIFF{--}results follow from comparing the
sequences (\ses) for \DIFF\ and \TOP\ using 
the Kirby\and Siebenmann action of $\br{X/\partial X}{\TOP/O}$
on both the normal maps and the structure sets
\diffSb\ needs an additional remark.
The outline above shows that a quotient of 
$H_0(B\pi_1;\cy2)$ acts freely on the 3{--}dimensional
structure set and this quotient can be compared with the
quotient for fundamental group with $\cy2$ and $w_1$
non{--}trivial.

\newsec{A Construction of Novikov, Cochran \and Habegger.}
\Remembersection\sectionfive

As we have seen above, the stable structure set 
in the simply connected case, while finite,
can be arbitrarily large. 
However, Freedman, \cite\bFa, says
that there are either one or two manifolds in each
homotopy type.
The resolution of this conundrum is the following.

Let $\hmty(X;\rel\ \partial X)$ denote
the group of degree one, simple homotopy automorphisms of
$X$, $\ell\colon (X,\partial X)\to (X,\partial X)$,
with $\ell\vert_{\partial X}=1_{\partial X}$.
Let $\ell$
act on $f\colon (M,L)\to(X,\partial X)\in\Scat(X;\rel\ h)$
via composition: 
$$\ell\autact f\colon (M,L)\RA{f} 
(X,\partial X)\RA{\ell} (X,\partial X)\ .$$
This group, $\hmty(X;\rel\ \partial X)$, 
acts on the stable structure sets, and even on each of the
$\wtScat{r}(X;\rel\ h)$, as follows.	
If $\ell\in\hmty(X;\rel\ \partial X)$,
there is a well{--}defined element in 
$\hmty(X\#rS^2\times S^2;\rel\ \partial X)$,
$\ell\#{id}\colon X\#r(S^2\times S^2)\to 
X\#r(S^2\times S^2)$
and we let $\ell$ act on $f\colon M\to X\#r(S^2\times S^2)$ in
$\wtScat{r}(X;\rel\ h)$ as the composite
$\ell\autact f\colon
M\RA{f} X\#r(S^2\times S^2)
\RA{\ell\#r{\rm id}}X\#r(S^2\times S^2)$.
The maps $\wtScat{r}(X;\rel\ h)\to
\hskip1pt\wtScat{r+1}(X;\rel\ h)$
and the maps from the \DIFF\ to the \TOP\ structure sets 
are equivariant with respect to these actions, 
so there are also actions on the stable
structure sets.

The set of \CAT{--}manifolds homotopy equivalent
to $X$, rel $h$, is just the orbit space
of this action.
The action preserves the stable triangulation obstruction,
so there is a set map
$$k\colon \Stop(X;\rel\ h)/\hmty(X;\rel\ \partial X)\to\cy2$$
and Freedman's classification follows from 
\cCRA\ below that $k$ is injective
in the simply connected case
plus the discussion of the image of $k$ in \nonE\ above.
Check that the embedding of $\hmty(X;\rel\ \partial X)$
in $\hmty(X\#S^2\times S^2;\rel\ \partial X)$ defined by
$\ell\mapsto \ell\#1_{S^2\times S^2}$ defines an action
of $\hmty(X;\rel\ \partial X)$ on $\barScat(X;\rel\ h)$.
\twotheorems{\thmNM}{\thmNMS}\ below give a partial calculation of
$\barScat(X;\rel\ h)/\hmty(X;\rel\ \partial X)$.

Let $X$ be a \CAT{--}manifold and use the identity
as a base point in $\Scat(X;\rel\ h)$.
Brumfiel, \cite\bGBr, shows that, in $[X/\partial X,\gcat]$, 
\Ndef\afone{}
$$\N_{\N_{1_X}}(\ell\autact f)=
\N_{\N_{1_X}}(\ell)+
(\ell^{-1})^\ast\bigl(\N_{\N_{1_X}}(f)\bigr)\ .\leqno(\afone)$$

\noindent
A similar formula holds for the action on the stable
structure sets.
Observe that any $\ell\in\hmty(X;\rel\ \partial X)$
preserves $w_1(X)$ and so induces an automorphism of
the Wall group $\LWG sn(\Z[\pi_1(X)],w_1(X))$.
One can check that with these definitions the
sequences (\ses) are $\hmty(X;\rel\ \partial X)$
equivariant.

There is a construction due to Novikov, \cite\bNo,
with the details finally worked out by 
Cochran \and Habegger, \cite\bCH.
Given any $\alpha\in \pi_2(X)$, let $\ell_\alpha$
denote the following composite
$$X\to X\vee S^4\RA{1_X\vee\eta^2}
X\vee S^2\RA{1_X\vee\alpha}X$$
where $\eta^2\in\pi_4(S^2)=\cy2$ denotes
the non{--}trivial element and the map
$X\to X\vee S^4$ just pinches the boundary
of a disk in the top cell to a point.

One point of Cochran \and Habegger's paper is to
compute the normal invariant of $\ell_\alpha$.
This result requires no fundamental group hypotheses 
and yields:

\Ndef\thmCH{Theorem}
\proclaim \thmCH.
$$\N_{1_X}(\ell_\alpha)=
(1+\langle v_2(X),\alpha\rangle)\bar\alpha$$
where 
$\bar\alpha\in\br{X/\partial X}{\gtop}$ 
denotes the image of $\alpha$
in $H_2(X;\cy2)\subset\br{X/\partial X}{\gtop}$ and 
$\langle v_2(X),\alpha\rangle\in\cy2$
denotes the evaluation of the cohomology
class on the homotopy class.

\remarks
Since $\ell_\alpha$ can be checked to induce the
identity on $[X/\partial X,\gcat]$, this formula 
and (\afone) determine the action of $\ell_\alpha$
on the \TOP{--}normal maps.
If $X$ is oriented, the \DIFF{--}normal maps
are a subset of the \TOP\ ones, so this formula
determines the action on the \DIFF{--}normal maps as well.
In the non{--}orientable case, there is a $\cy2$ in the
kernel of the map from the \DIFF{--}normal maps
to the \TOP\ ones and the Novikov{--}Cochran{--}Habegger
formula does not determine the normal invariant.

Let $\hmtyone(X; \rel\ \partial x)$
denote the subgroup of $\hmty(X;\rel\ \partial X)$
generated by the $\ell_\alpha$.

\Ndef\thmNM{Theorem}
\proclaim\thmNM.
$$\barStop(X;\rel\ h)/\hmtyone(X;\rel\ \partial X)
\RA{\N}\cases{
{\cal K}_4(\pi_1(X),w_1(X))&if $v_2(\widetilde X)=0$\cr
{\cal K}_4(\pi_1(X),w_1(X))\oplus\cy2&
if $v_2(\widetilde X)\neq0$\cr}
$$
is onto.
In the second case, the stable triangulation obstruction
is onto the $\cy2$: in the first case, $k$ may or may not
be constant as discussed in \nonE\ above.
Moreover
$\bar Q_5(\pi_1(X),w_1(X),D_X)$ acts transitively
on the orbits of this map.

\remark
\thmNM\ shows that except for a $\cy2$ related to stable
triangulation, there is an upper bound for
$\barStop(X;\rel\ h)/\hmtyone(X;\rel\ \partial X)$
which depends only on ``fundamental group data''.

\Ndef\cCRA{Corollary}
\proclaim \cCRA.
Suppose $\bar Q_5(\pi_1(X),w_1(X),D_X)=0$ and 
${\cal K}_4(\pi_1(X),w_1(X))=0$.
Then the set
$\barStop(X;\rel\ h)/\hmty(X;\rel\ \partial X)=
\barStop(X;\rel\ h)/\hmtyone(X;\rel\ \partial X)$
has one element 
if $\widetilde X$ is \spin,
and two elements with different
triangulation obstructions if it is not.
Any simple homotopy equivalence $f$
is homotopic to the composition of a homeomorphism
and an element in $\hmtyone$.

Notice that the action of $\bar\gamma$ on $\barSdiff$
preserves the $\hmtyone$ orbits, so
\Ndef\thmNMS{Theorem}
\proclaim \thmNMS.
The group $\bar\gamma$ acts on
$\barSdiff(X;\rel\ h)/\hmtyone(X;\rel\ \partial X)$
and the orbit space injects into
$\barStop(X;\rel\ h)/\hmtyone(X;\rel\ \partial X)$.

The action by the full group, $\hmty$, is more subtle and
often involves the homotopy of $X$, not just ``fundamental group data''.
Let $X$ be a \TOP{--}manifold and define
$$\hmtytwo(X;\rel\ \partial X)=\{\ell\in\hmty(X;\rel\ \partial X)
\ \vert\ \N_{\N(1_X)}(\ell)=0 \ and\ \ell_\ast=1_{\pi_1(X)} \}\ .$$
It follows from Brumfiel's formula (\afone) that
$\hmtytwo$ is a subgroup of $\hmty$.
Let $\hmtybig$ denote the subgroup of $\hmty$
generated by $\hmtyone$ and $\hmtytwo$.
\twotheorems{\thmNM}{\thmNMS}\ 
continue to hold with $\hmtybig$
replacing $\hmtyone$.
The actual homotopy type of $X$ can be seen to effect
$\barScat(X;\rel\ h)/\hmtybig(X;\rel\ \partial)$ via the following
observation.
The evident map
$\barScat(X;\rel\ h)\to\barScat(X\#S^2\times S^2;\rel\ h)$
induces a map
$$\iota_X\colon\barScat(X;\rel\ h)/\hmtybig(X;\rel\ \partial)\to
\barScat(X\#S^2\times S^2;\rel\ h)/
\hmtybig(X\#S^2\times S^2;\rel\ \partial)\ .$$

\noindent
Let $\WSEcat(X;\rel\ h)$ denote the limit of the maps
$\iota_X, \iota_{X\#S^2\times S^2}, \cdots,
\iota_{X\#rS^2\times S^2},\cdots$.

\Ndef\bTHMa{Theorem}
\proclaim \bTHMa.
The evident quotient of the normal map,
$$\WSEtop(X;\rel\ h)
\RA{\N}\cases{
{\cal K}_4(\pi_1(X),w_1(X))&if $v_2(\widetilde X)=0$\cr
{\cal K}_4(\pi_1(X),w_1(X))\oplus\cy2&
if $v_2(\widetilde X)\neq0$\cr}$$
is a bijection.
If $k^{-1}(0)\neq\emptyset$, then
$\WSEdiff(X;\rel\ h)\to k^{-1}(0)$ is a bijection.

\remarks
The stable triangulation obstruction
is onto the $\cy2$ if $v_2(\widetilde X)\neq0$:
otherwise, $k$ may or may not
be constant as discussed in \nonE\ above.
Note that ${\cal K}_4$ is always a $\cy2$ vector
space of dimension at most $H_2(B\pi;\cy2)\oplus\cy2$,
and hence finite.
If $\hat Q_5(\pi_1(X), w_1(X))$ is finitely generated,
then there exists an $r$ such that
$$\barScat(X\#rS^2\times S^2;\rel\ h)/
\hmtybig(X\#rS^2\times S^2;\rel\ \partial X)\to\WSEcat(X;\rel\ h)$$
is a bijection.

\newsec{Examples.}
\Remembersection\sectionsix

Here are some calculations for some specific manifolds.
The quoted values of 
$\bar Q_5$, ${\cal K}_4$ and $\bar\gamma$
can be obtained from \tabo, after noting 
\QRt\ applies so $\bar Q_5=\hat Q_5$ and $\bar\gamma=\hat\gamma$.

\def\example#1.{\par\noindent{\bf Example: #1.} }
\example $RP^4$.
Here $\pi=\cy2$ and $w_1$ is an isomorphism.
Then
$\bar Q_5=0$, ${\cal K}_4=\cy2$ and $\bar\gamma=\cy2$.
For the normal maps, 
$\br{RP^4}{\gtop}=\cy2\oplus\cy2$; $\br{RP^4}{\gdiff}=\cy4$
and it is a useful exercise to understand how the \DIFF\ and \TOP\
versions of sequence (\ses) work in this case without
relying on the general theory.

Hence $\barStop(RP^4)=\cy2$ and $k$ is a bijection.
The non{--}triangulable example was constructed by Ruberman, \cite\bRu,
using only Freedman's simply connected results.
In the smooth case, $\barSdiff(RP^4)=\cy2$ as well, but the map
from the smooth to the topological sets takes both elements
of the smooth set to one element in the topological set.
Cappell\and Shaneson, \cite\bCSfh,
constructed an element in $\Sdiff(RP^4)$
which hits the ``other element'' in $\barSdiff(RP^4)$.

\example
$S^3\times S^1$.
Here $\pi=\Z$ and $w_1$ trivial.
Then
$\bar Q_5=0$, ${\cal K}_4=0$ and $\bar\gamma=\cy2$.

It follows that $\barStop(S^3\times S^1)$ is one point
and $\barSdiff(S^3\times S^1)$ is two points.
The ``other element'' in $\barSdiff(S^3\times S^1)$ was
constructed in $\wtScat 1(S^3\times S^1)$ by 
Scharlemann, \cite\bMS.
It is an open question as to whether this element is
in the image from $\Sdiff(S^3\times S^1)$.

\example
$S^3\widetilde\times S^1$.
Here $\pi=\Z$ and $w_1$ non{--}trivial.
Then
${\cal K}_4=0$, $\bar Q_5=0$ and $\bar\gamma=\cy2$.

Hence $\barStop(S^3\widetilde\times S^1)$ consists of
one point, while $\barSdiff(S^3\widetilde\times S^1)$
consists of two points, distinguished by the smooth
normal invariant.
In this case, Akbulut, \cite\bAK, constructed the ``other element''
in $\wtSdiff 1(S^3\widetilde\times S^1)$.

\remark
If one could find a manifold to show
$\psi_{\DIFF}$ were onto for $S^3\widetilde\times S^1$,
Lashof\and Taylor, \cite\bLT, observed that 
$\bar\gamma$ would act freely on $\Sdiff(X;\rel\ h)$
as soon as this structure set is non{--}empty.
It does act freely on $\Sdiff(X\#S^2\times S^2;\rel\ h)$.
Cappell\and Shaneson's work, \cite\bCSfh, shows that $\bar\gamma$
acts freely on $\Sdiff(X;\rel\ h)$ if $\pi_1(X)=\cy2$ 
and $w_1$ is non{--}trivial.

\example
$RP^3\times S^1$.
Here $\pi=\cy2\times\Z$  and $w_1$ is trivial.
Then
${\cal K}_4=\cy2\oplus\cy2$,
$\bar Q_5=\Z$ and $\bar\gamma=\cy2$.

The manifold $RP^3\times S^1$ is \spin, so
$\barSdiff(RP^3\times S^1)\to 
\barStop(RP^3\times S^1)$ is onto.
There are two elements in $\barSdiff(RP^3\times S^1)$
over each element of $\barStop(RP^3\times S^1)$.
Each orbit of the Wall group has countable many elements
falling into 4 orbits, distinguished by the normal invariant.
For some $r$,
$\barScat(RP^3\times S^1\# rS^2\times S^2)/
\hmty(RP^3\times S^1\# rS^2\times S^2)$
contains at most 4 elements.

\remembersection\Sfive
\newsec{The Topological Case in General.}
\Remembersection\sectionseven

In a series of papers, Freedman,
\cite\bFa, \cite\bFb, \cite\bFQ,
showed that the high dimensional theory of surgery
and the high dimensional $s${--}cobordism theorem
hold in the \TOP{--}category in dimension 4 for 
certain fundamental groups.
As of this writing, there are no known failures
of either surgery theory or the $s${--}cobordism theorem
in the \TOP{--}category in dimension 4.

We say {\sl \CAT{--}surgery works in dimension $n$ for 
fundamental group $\pi$},
provided that, for any $n${--}dimensional Poincar\'e space 
$X$ with fundamental group $\pi$, the map
$$\psi_{\CAT}\colon\Scat(X;\rel\ h)\to\barScat(X;\rel\ h)$$
is a surjection;
we say the {\sl{\CAT}{--}$s${--}cobordism works in
dimension $n$ for fundamental group $\pi$},
provided that, for any $n${--}dimensional 
Poincar\'e space $X$ 
with fundamental group $\pi$, the map
$$\psi_{\CAT}\colon\Scat(X;\rel\ h)\to\barScat(X;\rel\ h)$$
is an injection.

The first of Freedman's theorems is

\proclaim Theorem.
\TOP{--}surgery and the \TOP{--}$s${--}cobordism theorem 
work in dimension $4$
for trivial fundamental group.

It took some work to get to this statement.
Freedman began with the simply connected, smooth case,
building on work of Casson, \cite\bAC.
By showing that Casson handles were topologically standard,
Freedman showed that surgery theory and the h{--}cobordism
theorem held {\sl topologically\/} for simply connected,
{\sl smooth\/} manifolds.

Quinn melded these results with his controlled results
to prove 
$\pi_i(\TOP(4)/O(4))=0$, $i=0, 1, 2$; $i=0$ is the
annulus conjecture in dimension 4.
Lashof\and Taylor, \cite\bLT, showed 
$\pi_3(\TOP(4)/O(4))=\cy2$ 
and reproved Quinn's result for $i=2$.
Finally, Quinn, \cite\bQa, showed $\pi_4(\TOP(4)/O(4))=0$, 
thus computing the last of the 
``geometrically interesting'' homotopy groups.
Using these results, Quinn, \cite\bFQ, then went on to show
that transversality works inside of topological 
4{--}manifolds.
Freedman had already completed a program of
Scharlemann, \cite\bMSa, by showing
that transversality worked in other dimensions
when the expected dimension of the result was 4.
After this, the standard geometric tools were available 
in dimension 4
and \TOP{--}surgery and the \TOP{--}s{--}cobordism 
theorem now worked for trivial fundamental group.

Freedman, \cite\bFb, then introduced capped{--}grope 
theory which he used to extend the fundamental groups
for which \TOP{--}surgery theory and 
the \TOP{--}$s${--}cobordism theorem work.
There is a nice general result, explained in \cite\bFT.
Following that exposition, we say that a group, 
$\pi$, is NDL, 
for {\sl Null Disk Lemma}, provided that, 
for any height 2 capped grope, $G$,
and any homomorphism, $\psi\colon \pi_1(G)\to\pi$,
we can find an immersed core disk, so that all the
double point loops map to $0$ under $\psi$.

\Ndef\thmNDL{Theorem}
\proclaim\thmNDL.
If $\pi$ is an NDL group, then \TOP{--}surgery and the 
\TOP{--}$s${--}cobordism theorem
work in dimension $4$ for $\pi$.

Freedman\and Teichner, \cite\bFT, check that any extension
of an NDL group by another NDL group is itself NDL,
and they check  that a direct limit of NDL groups is NDL.
Transparently, subgroups of NDL groups are NDL, and, 
since $\pi_1(G)$ is a free group,
quotients of NDL groups are NDL.
Hence subquotients of NDL groups are NDL and a group
is NDL iff all its finitely{--}generated subgroups are.
Finally, the main result of \cite\bFT, is
\Ndef\sexp{Theorem}
\proclaim \sexp.
Groups of subexponential growth are NDL.

It is possible that all groups are NDL.
Since any finitely{--}generated group is a subquotient
of the free group on 2 generators, all groups are NDL iff
the free group on 2 generators is.
An equivalent formulation, which might make the result
seem less likely, is that all groups are NDL iff
each height 2 capped grope has an immersed core disk with
all double point loops null homotopic.

Among the groups satisfying NDL are the finite groups, $\Z$,
$\Q$ and nilpotent groups.
There do exist nilpotent groups of exponential growth 
\ecite{\bPL}{Problem 4.6}.

Free groups on more than one generator are not known to be
NDL and this causes a great many other 
geometrically interesting groups to be on the unknown list.
Surface groups for genus 2 or more are examples of
such groups.
The free product of two groups, neither of which is 
trivial, is either
$\cy2\ast\cy2=\Z\>\semidirect\cy2$ (and is NDL)
or has a free subgroup of rank 2
(and is not known to be NDL).
Hence the fundamental groups of most connected sums
of 3{--}manifolds  are not known to be NDL.
Among the irreducible 3{--}manifolds, 
many are hyperbolic by Thurston, \cite\bTh,
and many of these have incompressible surfaces:
the fundamental groups of such manifolds are 
not known to be NDL.
Even if some group fails to be NDL, 
it is not clear that \TOP{--}surgery
must therefore fail for it.
In Freedman \& Quinn \cite\bFQ\ there is a
different condition whose truth would yield surgery and the
$s${--}cobordism theorem. 
It is possible that this condition could yield results even if
the Null Disk Lemma were to fail.
Quinn \cite{\bQz}\ has a nice discussion of the current state of
affairs regarding the groups for which 
surgery and/or the $s${--}cobordism theorem works.

In dimension 3, \TOP{--}surgery sometimes holds for trivial reasons:
for fundamental group trivial, $\Z$ (with either $w_1$)
or groups satisfying the Borel conjecture 
(which is conjectured to hold for all irreducible 
3{--}manifold groups), the stable \TOP\ structure set is trivial
and so \TOP{--}surgery holds.
For non{--}trivial, finite fundamental group,
\TOP{--}surgery fails for closed manifolds.
As an example $\barStop(RP^3)=\Z$ but 
Casson, \cite\bAM, shows that 
$\psi_{\TOP}\colon\Stop(RP^3)\to\barStop(RP^3)$
cannot hit an element of odd order
since the double cover of any such
element would be a homotopy 3{--}sphere of
Rochlin invariant 1.
This line of argument works for any other finite fundamental group.
The \DIFF{--}case is even worse 
since $\barSdiff(S^3)=\cy2$ and
the result of Casson's used above also shows that 
$\psi_{\DIFF}$ is not onto.

To say that the $s${--}cobordism theorem holds
in dimension 3 is a bit of a misnomer.
If \TOP{--}surgery works for $\pi_1(X)$, then
two elements in $\Stop(X;\rel\ h)$ which hit 
the same element in 
$\barStop(X;\rel\ h)$ differ by an $s${--}cobordism.
However, as we saw above, we do not know whether 
\TOP{--}surgery holds for many 3{--}manifold groups
and hence we do not usually know 
that there is an $s${--}cobordism between the two elements.
Still, we retain the terminology despite its drawbacks.

For $\pi$ trivial, the $s${--}cobordism theorem
holds in dimension 3 iff the Poincar\'e conjecture holds.
In general, the $s${--}cobordism fails
in the strict sense that there are 4{--}dimensional
\TOP{--}s{--}cobordisms which are not products.
The first such examples are due to 
Cappell\and Shaneson, \cite\bCSh, with a much
larger collection of examples worked out by 
Kwasik\and Schultz, \cite\bKSh.
Surprisingly, there are no counterexamples known 
to us of the smooth $s${--}cobordism theorem failing
in dimension 4, but this is probably due to our
inability to construct smooth $s${--}cobordisms.

It may be worth remarking that two 4{--}dimensional
results from the past now can be pushed down
one dimension.
Barden's old observation that an $h${--}cobordism
from $S^4$ to itself is a smooth product can be
made again to observe any $h${--}cobordism from
$S^3$ to itself is a topological product.
Thomas's techniques, \cite\bCTa, can be applied
to show that any 4{--}dimensional
$s${--}cobordism with NDL fundamental
group is invertible.

There has been a great deal of work using Freedman's
ideas to attack old problems in four manifolds.
A complete survey of such results 
would require more than our allotted space.
Here are some examples which have lead to further work.
Hambleton, Kreck\and Teichner classify
non{--}orientable 4{--}manifolds with fundamental
group $\cy{2}$, \cite\bHKT.
Hambleton\and Kreck also classify orientable 
4{--}manifolds with fundamental group $\cy{N}$,
\cite\bHK, as the start of a general program
to extend Freedman's simply connected classification
to manifolds with finite fundamental group.
Kreck's reformulation of surgery theory works very well here,
\cite\bMKy.

Lee\and Wilczynski, \cite\bLW, have largely solved the problem
of finding a minimal genus surface representing a
2{--}dimensional homology class in a
simply connected 4-manifold.
Askitas \cite{\bAsk}\ and \cite{\bAskb}\ considers some cases 
of trying to represent several homology classes at once.

The slicing of knots and links is an active area as well.  The first
results here were negative.  Casson \& Gordon's examples \cite{\bCG}\
\cite{\bCGa}\ of algebraically slice knots that were not slice showed
that there does not exist enough embedding theory in dimension four to
do $\Gamma${--}group surgery in the style of Cappell \& Shaneson \cite{\bCSz}.

One of Freedman's striking results \cite{\bFb}\ is that knots of Alexander
polynomial $1$ are topologically slice.  
Casson \& Freedman \cite{\bCF}\ found links which would be slice if and only
if surgery theory worked in dimension $4$ for all groups.

On the other hand, it was known in the 1970's to Casson (and others?)
that in a smooth 4-manifold $M$ with no 1-handles, the only
obstruction to representing a characteristic class of square one by a
PL embedded 2-sphere with one singularity with link a knot of
Alexander polynomial one, was the Arf invariant of the knot (that is,
$\sigma M \equiv 1 \bmod{16} $).  
Once Donaldson showed that non-diagonal definite forms were not
realized by smooth 4-manifolds, then in $\CP^2$ blown up at 16 points,
any characteristic class of square 1 cannot be represented by a
smoothly embedded 2-sphere. 
Hence there must be an Alexander polynomial
one knot which is not smoothly slice in a homology 4-ball. (See
Problem 1.37, page 61 in  \cite{\bPL})

\remembersection\Ssix
\newsec{The Smooth Case in Dimension 4.}
\Remembersection\sectioneight
Shortly after Freedman's breakthrough in 1981,
Donaldson made
spectacular progress in the smooth case.
We soon learned that neither \DIFF{--}surgery nor the
\DIFF{--}$s${--}cobordism theorem holds, even for simply connected
smooth manifolds.
In the next fifteen years, we learned a great deal more,
but the overall situation has only become more complex
from the point of view of surgery theory.

\vskip4pt\noindent{\bf Existence:}
Donaldson's first big theorem, \cite\bSDa,
severely limited the forms which
could be the intersection form of a 
smooth, simply connected 4-manifold.
Any form can be stably realized and as soon as the form is
indefinite, they are completely classified.
In the \spin\ case, the forms are $2mE_8\oplus rH_2$
where $E_8$ is the famous definite even form of signature $8$
and $H_2$ is the dimension $2$ hyperbolic.
Donaldson, \cite\bSDb, proved that if $m=1$, then $r\geq3$, 
and there is a conjecture,
the $11/8${--}th's conjecture
($b_2/\vert\sigma\vert\geq 11/8$), 
which says that $r\geq3m$ in general.
At this time Furuta, \cite\bFru,
has proved the $10/8${--}th's conjecture,
which says that $r\geq2m$.
See \cite\bPL, Problems 4.92 and 4.93.
In particular, there exists a simply connected,
\TOP\ manifold, $M_{2mE_8}$ with form $2mE_8$;
from \diffS, $\barSdiff(M_{2mE_8})=16m(\cy2)$, but
$\wtSdiff{r}(M_{2mE_8})=\emptyset$ for $r<2m$.
In the simply connected case, we also know that, 
for each integer $r\geq 0$, either
$\wtSdiff{r}(M)=\emptyset$ or else 
$\psi^r_{\DIFF}\colon \wtSdiff{r}(M)
\to\barSdiff(M)$ is onto.

Scharlemann, \cite\bMS, showed 
$\psi^1_{\DIFF}\colon \wtSdiff{1}(S^3\times S^1)\to
\barSdiff(S^3\times S^1)=\cy2$
is onto: $\Sdiff(S^3\times S^1)$ is certainly
non{--}empty, but as of
this writing, $\psi_{\DIFF}$ is not known to be onto.
Wall, \ecite\bCTC{\S 16}, shows all homotopy
equivalences are homotopic to diffeomorphisms,
so $\hmty(S^3\times S^1)$ acts trivially 
on the smooth structure set.
Interestingly, a folk result of R.~Lee,
\cite\bCSs, says that
$\hmty(S^3\times S^1\# S^2\times S^2)$ acts transitively on
$\wtSdiff{1}(S^3\times S^1)$.

The above gives many examples of simply connected smooth manifolds
which topologically decompose as connected sums, but have no
corresponding smooth decomposition.
Works of Freedman\and Taylor, \cite\bFTx,
and Stong, \cite \bRS, show that one can still mimic
this decomposition by decomposing along homology 3{--}spheres into
simply connected pieces.

\vskip4pt\noindent{\bf Uniqueness:}
Donaldson, \cite\bSDc, also proved that 
the $h${--}cobordism theorem fails for
smooth, $5${--}dimensional, simply connected 
$h${--}cobordisms.
Note however that a smooth h-cobordism between simply connected
4-manifolds is unique up to diffeomorphism, \cite\bMKx.
There is another classification theorem of simply connected $h${--}cobordism
due to Curtis, Freedman, Hsiang\and Stong, \cite\bCFHS,
in terms of Akbulut's corks, \cite\bKC, \cite\bACork, \cite\bAstudent.

We know of no case in which $\psi_{\DIFF}$ is not 
$\infty${--}to{--}one
and we know of no case where all the elements in 
$\Sdiff(M)$ have been described.
The smooth Poincar\'e conjecture, unresolved at the time
of this writing, says $\Sdiff(S^4)$ has one element.
The uniqueness result for $\R^4$ is known to fail 
spectacularly, \cite\bBG, \cite\bDF.
In contrast to the existence question,
where we know examples for which we need arbitrarily
many $S^2\times S^2$'s before
a particular stable element exists, for all we know,
$\wtSdiff{r}(M)\to\barSdiff(M)$ and 
$\wtSdiff{r}(M)\to\wtSdiff{r+1}(M)$
have the same image.
Some works of Mandelbaum\and Moishezon, \cite\bMM, and Gompf, \cite\bBGx,
give many examples in which this one{--}fold stabilization suffices.

It follows from Cochran\and Habegger, \cite\bCH,
that the group of homotopy automorphisms of a closed, simply connected
4{--}manifold, $M$, is the semidirect product of the Novikov maps,
$\hmtyone(M)$, and the automorphisms of $H_2(M;\Z)$ which preserve
the intersection form.
Moreover, Cochran\and Habegger show that all the non{--}trivial elements
of $\hmtyone(M)$ are detected by normal invariants and so are
not homotopic to homeomorphisms.
Now it follows, as observed by Freedman, \cite\bFa, that 
the automorphisms of $H_2(M;\Z)$ which preserve the intersection form
are realized by homeomorphisms, unique up to homotopy.
Further work by Quinn, \cite\bQa, shows that they are in fact unique
up to isotopy.

When $M$ is also smooth and of the form $P\# S^2\times
S^2$, Wall, \cite\bWd, and Freedman\and Quinn, \cite\bFQ, showed that
any homeomorphism is isotopic to a diffeomorphism.
But when $M$ is
not of the form $P\# S^2\times S^2$, then there are often severe
restrictions on realizing a homotopy equivalence by a diffeomorphism
due to the existence of basic classes in $H^2 (M, \Z)$.
These classes
were defined for Donaldson theory by Kronheimer\and Mrowka,
\cite\bKMa, \cite\bKMb.
Conjecturally equivalent basic classes were
also defined using Seiberg-Witten invariants, \cite\bSW, and these
classes were shown to be equivalent by Taubes, \cite\bTa, to classes
defined via Gromov's pseudoholomorphic curves.
Although the set of
basic classes can be as simple as the zero class in $H^2 (M,\Z)$ for
the $K3$ surface, the classes can be as complicated as Alexander
polynomials are, \cite\bFS.
The isometry induced on
$H^2(M;\Z)$ by a diffeomorphism must take each basic class to 
$\pm($ a, possibly different, basic class$)$.

There can be further restrictions, beyond
those determined by the basic classes, to realizing homotopy
equivalences by diffeomorphisms.
For example, any $K3$ surface has additional restrictions, see
\ecite{\bDK}{Corollary 9.14, p.345}.
The homeomorphism of $K3$ which is the identity except on 
an $S^2\times S^2$
summand and is $antipodal\times antipodal$ on the 
$S^2\times S^2$ summand
cannot be realized by a diffeomorphism.
However, it follows from \cite\bFM, that a subgroup of finite index in 
the group of isometries of the intersection form
of $K3$ is realized by 
diffeomorphisms.

As of this writing, work in the smooth case is 
continuing at a feverish pace and is hardly ripe for a survey.
For many smooth manifolds we now know the minimal genus
smooth embeddings representing any homology class; see Kuga
\cite{\bKuga}, Li \& Li \cite{\bLL}\cite{\bLLb}, Kronheimer \& Mrowka \cite{\bKMa},
and Morgan, Szab\'o \& Taubes \cite{\bMST}.
Some work on simultaneous representation of several classes
in the smooth case is in \cite{\bAsk}.
The xxx Mathematics Archive at Los{--}Alamos (see
http://front.math.ucdavis.edu/ ) is a useful resource for those
wishing to remain current.

\global\truesection=200

\centerline{\bf Bibliography}
\vskip10pt
\showbiblist

\end